\documentclass[10pt]{article}
\usepackage[latin1]{inputenc}
\usepackage{amsmath}
\usepackage{amsfonts}
\usepackage{bm}
\usepackage{amssymb}
\usepackage{graphicx}
\usepackage{subcaption}
\usepackage{multirow}
\usepackage{longtable}
\usepackage{authblk}
\usepackage{rotating}
\usepackage{lscape}
\usepackage{fancyvrb}
\usepackage{tabu}
\usepackage[hidelinks]{hyperref}
\usepackage{xcolor}
\hypersetup{
    colorlinks,
    linkcolor={black!50!black},
    citecolor={black!50!black},
    urlcolor={blue!80!black}
}

\allowdisplaybreaks
\newcommand{\g}[1]{\Gamma(#1)}
\newcommand{\p}[2]{(#1)_{#2}}

\date{}
\newcommand{\linefill}{
  {-}\mkern-7mu
  \cleaders\hbox{$\mkern-2mu-\mkern-2mu$}\hfill
  \mkern-7mu{-}%
}

\begin{document}

\begin{titlepage}

\begin{center}

\begin{flushright}

\end{flushright}

\begin{flushright}

\end{flushright}

\begin{flushright}

\end{flushright}

{\Large\bf On the evaluation of the \\} 

\medskip

{\Large\bf Appell $F_2$ double  hypergeometric function\\}

\vspace{1.5cm}

{\bf B. Ananthanarayan$^{a\ast}$, Souvik Bera$^{a\dagger}$, S. Friot$^{b,c\ddagger}$, O. Marichev$^{d\diamond}$ and Tanay Pathak$^{a\star}$}\\[1.5cm]
{$^a$ Centre for High Energy Physics, Indian Institute of Science, \\
Bangalore-560012, Karnataka, India}\\[0.5cm]
{$^b$ Universit\'e Paris-Saclay, CNRS/IN2P3, IJCLab, 91405 Orsay, France } \\[0.5cm]
{$^c$ Univ Lyon, Univ Claude Bernard Lyon 1, CNRS/IN2P3, \\
 IP2I Lyon, UMR 5822, F-69622, Villeurbanne, France}
\\[0.5cm]
{$^d$ Wolfram Research, Inc., 100 Trade Center Drive, Champaign, IL 61820-7237, USA}  \\[2.cm]
\end{center}

\begin{abstract}
The transformation theory of the Appell $F_2(a,b_1,b_2;c_1,c_2;x,y)$ double hypergeometric function is used to obtain a set of series representations of $F_2$ which provide an efficient way to evaluate $F_2$ for real values of its arguments $x$ and $y$ and generic complex values of its parameters $a,b_1, b_2, c_1$ and $c_2$ (i.e. in the nonlogarithmic case). This study rests on a classical approach where the usual double series representation of $F_2$ and other double hypergeometric series that appear in the intermediate steps of the calculations are written as infinite sums of one variable hypergeometric series, such  as the Gauss $_2F_1$ or the $_3F_2$, various linear transformations of the latter being then applied to derive known and new formulas. Using the three well-known Euler transformations of $F_2$ on these results allows us to obtain a total of 44 series which form the basis of the \textit{Mathematica} package \textsc{AppellF2}, dedicated to the evaluation of $F_2$. A brief description of the package and of the numerical analysis that we have performed to test it are also presented.

\end{abstract}

\vspace{3cm}
\small{$\ast$ anant@iisc.ac.in }

\small{$\dagger$ souvikbera@iisc.ac.in}

\small{$\ddagger$ samuel.friot@universite-paris-saclay.fr}

\small{$\diamond$ oleg@wolfram.com}

\small{$\star$ tanaypathak@iisc.ac.in}

\end{titlepage}

\section{Introduction}

The study of multiple hypergeometric functions, which appear in many domains of the physical and mathematical sciences, mainly began in 1880 with Appell who introduced the famous four double hypergeometric functions $F_1, F_2, F_3$ and $F_4$ that carry his name and are generalizations of the Gauss hypergeometric $_2F_1$ function. Since then, a huge development of this field, by the extensive study of many classes of multiple hypergeometric functions, led it to become a classical branch of mathematics. On the other hand, the implementation of the automatic evaluation of these multivariable functions in softwares dedicated to mathematics is a field of investigations which is nearly virgin. For instance, in \textit{Mathematica} \cite{Mathematica} only the Appell $F_1$ has been coded, and in \textit{Maple} \cite{Maple}, the four Appell functions are the only hypergeometric functions of more than one variable that are in-built functions (since 2017). 

One can point out several difficulties that may be at the origin of such a lack. One of them is that the integral representations of multiple hypergeometric functions are not always known and, when known, they are in general not valid for all values of the parameters of the hypergeometric functions that they represent. They can also be hard to compute numerically. 

An alternative way to handle multiple hypergeometric functions is to consider their series representations and, using transformation theory \cite{Bateman, Srivastava}, to obtain other series representations converging in other regions of the space of their variables, giving thereby analytic continuations of the starting point series. 
One interest of this approach is that the convergence properties of multiple hypergeometric series are independent of the values of their parameters (exceptional values of the parameters being excluded). Thus, one can use these series representations for numerical purpose when the Euler integral representations (or other integral representations) are not defined, or are unknown. Another advantage is that it is often easier to numerically compute series than integrals. However, one has to point out that beyond the case of double series, the convergence regions of multivariable hypergeometric series can be difficult to obtain. Moreover, to our knowledge, there is no systematic approach to derive transformations of these series that can collectively provide an evaluation of the corresponding multiple hypergeometric functions for all the possible values of their variables. 

In what concerns the analytic continuation of multiple hypergeometric series, a recent and important progress can be mentioned. In \cite{ABFGgeneral}, two of the authors of the present paper have developed, with other collaborators, a very efficient and systematic approach to analytically compute multiple Mellin-Barnes (MB) integrals. It is well-known that MB integrals are intimately linked to hypergeometric functions \cite{Appell, Exton, Marichev}. Multiple MB integrals are in fact one of the possible starting points for the study of hypergeometric functions of several variables \cite{TZ}. Therefore, by the study of appropriate classes of multiple MB integrals, the method of \cite{ABFGgeneral} opens promising horizons in hypergeometric functions theory and, in particular, for the determination of the analytic continuations of many classes of hypergeometric series, whatever the number of their variables is, in terms of others multivariable hypergeometric series. 
Obviously, a large number of applications can follow in physics, as already shown in the recent works \cite{Ananthanarayan:2020ncn} and \cite{Ananthanarayan:2020xpd} in the context of the study of Feynman integrals in quantum field theory.

Although many new results can be obtained from the powerful method developed in \cite{ABFGgeneral}, it cannot, alone, fully solve the difficult problem of finding the relevant set of transformations of a multiple hypergeometric function that will allow to evaluate it for all values of its variables. Indeed, we have mentioned above some possible difficulties in the derivation of the convergence regions of the new series representations obtained from transformation theory. But another problem can be that some hypergeometric functions of several variables do not have an obvious MB representation. Moreover, even if the latter can be obtained, 
 the evaluation of the MB integral in general shows that ``white regions" (as called in \cite{Ananthanarayan:2020acj}) appear in the multivariable space, where none of the obtained analytic continuations converges. Although some manipulations of the MB integral can lead to transformations and, thus, to other formulas (in addition to those obtained by a direct application of \cite{ABFGgeneral}), it is not clear that a systematic approach for these manipulations can be found. Therefore, in order to fully solve the problem of evaluating multivariable hypergeometric functions for all values of their variables, it may be necessary to rely on alternative approaches, in order to complete the results obtained from the MB approach.

A well-known example of such a situation in the context of quantum field theory involves the triple hypergeometric function of Lauricella $F_C$ type \cite{Berends:1993ee}. This particular function, which is the natural extension of the Appell $F_4$ double hypergeometric function, appears when one computes the two-loop sunset Feynman integral with four mass scales. It is easy to conclude from \cite{Berends:1993ee} that the analytic continuations of the $F_C$ triple series, derived from the Mellin-Barnes representation of the $F_C$ function, give access to a restricted region of its three variables space. Therefore, in order to obtain analytic expressions for the sunset outside this region, some transformations of the $F_C$ Lauricella series have been obtained in \cite{Ananthanarayan:2019icl}, using an alternative method. This method, which uses quadratic transformations of the Gauss $_2F_1$ hypergeometric series as intermediate steps in the derivation of new series representations for $F_4$ \cite{Ananthanarayan:2020xut} (and, as a by product, for $F_C$) can be seen as an extension of a classical work of Olsson \cite{Olsson-64} which focused on the question of the analytic contination of the Appell $F_1$ series and of its $F_D$ multivariable generalisation, using linear transformations of $_2F_1$. The approach of \cite{Ananthanarayan:2019icl,Ananthanarayan:2020xut} can however not give the full answer to the problem of finding series representations that can be used to evaluate the $F_C$ function for all values of its variables. 

Our aim in the present work is to explore Olsson's approach more systematically, taking the simpler case of the Appell $F_2$ double hypergeometric function as a theoretical laboratory, having in mind, among others, to come back at a later stage to the case of the Lauricella $F_C$ and other multivariable hypergeometric functions. 

The Appell $F_2$ double hypergeometric function is not an arbitrary choice, it has indeed a particular place in the set of the 14 complete double hypergeometric functions of order 2, which consist of the four Appell functions $F_i, (i=1,...,4)$ and the ten Horn functions $G_i, (i=1,..,3)$ and $H_i, (i=1,...,7)$. Indeed, it has been noticed in \cite{Erdelyi-1948} that, with the exception of $F_4, H_1$ and $H_5$, the Appell $F_2$ function can be related to any of the other Horn and Appell functions. These links can be obtained from the transformation theory of $F_2$ and are summarized in Chapter 5 of \cite{Bateman}.  
Hence, in the present work, by studying the linear transformations\footnote{Some work has been performed  by Olsson long ago on the study of the partial differential equations system of $F_2$ \cite{Olsson-f2}, whose solutions have been exhibited. However, the transformation formulas needed for the present work have not been derived in this reference.} of $F_2$ and by building the \textsc{AppellF2} \textit{Mathematica} package based on the obtained formulas and dedicated to its numerical evaluation, we provide the basis of a future \textit{Mathematica} package for the evaluation of all the double series above \cite{workinprogress}, with the exception of $F_4$, $H_1$ and $H_5$. These three lacking series will be considered separately in subsequent publications.

The plan of the paper is as follows. In Section \ref{F2_Sec}, we briefly list some of the well-known properties of the Appell $F_2$ function. In Section \ref{MB_F2Sec}, we perform a first analytic continuation study of the Appell $F_2$ series from the Mellin-Barnes approach \cite{ABFGgeneral}, which is completed in Section \ref{Olsson_Sec} following Olsson's method. This analysis, which yields 11 series representations of $F_2$, can be extended with the use of the three Euler transformations of $F_2$, allowing us to obtain a total set of 43 linear transformations of $F_2$, out of which 17 are needed, in addition to the usual series definition of $F_2$, to cover the $(x,y)$ space of the $F_2$ variables for real values of the latter, with the exception of a few points. This subset of 18 series representations of $F_2$ are recapitulated in the appendix along with the figures showing their corresponding regions of convergence for real values of the arguments. The mathematical expressions and convergence regions of the remaining 26 series can be obtained from our \textsc{AppellF2} package. These additional series increase the efficiency of the package  from the convergence perspective by enlarging the possible ways to compute $F_2$. 
Section \ref{package} is dedicated to the description of the \textsc{AppellF2} package 
 and its numerical tests, which are followed by the conclusions and the appendix.

\section{The Appell $F_2$ function \label{F2_Sec}}

\begin{figure}[h]
\centering
\includegraphics[width=6cm, height=6cm]{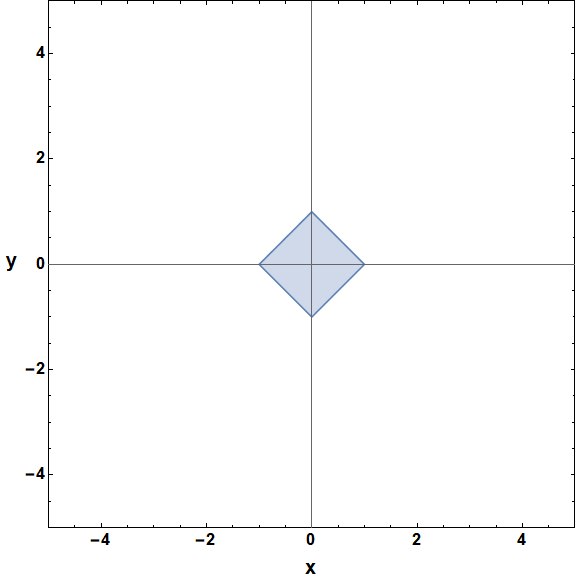}
\includegraphics[width=6cm, height=6cm]{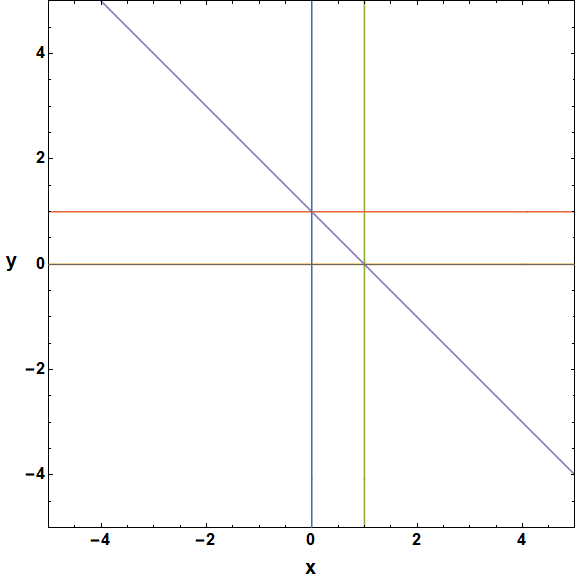}
\caption{\textit{Left}: Region of convergence (ROC) of the Appell $F_2$ series for real values of its arguments $x$ and $y$.\label{RocF2} \textit{Right}: Singular curves of Eq.(\ref{PDE}).}
\end{figure}

The Appell $ F_2$ double hypergeometric series is defined as \cite{Appell}
\begin{equation}\label{F_2original}
F_2(a,b_1,b_2;c_1,c_2;x,y)= \sum_{m=0}^{\infty}\sum_{n=0}^{\infty}\frac{  \p{a}{m+n}\p{b_1}{m} \p{b_2}{n}}{ \p{c_1}{m} \p{c_2}{n}m! n!}x^m y^n
\end{equation}
where $(a)_m=\frac{\Gamma(a+m)}{\Gamma(a)}$ is the Pochhammer symbol. The series in the RHS of Eq.(\ref{F_2original}) converges for $|x|+|y|<1$ 
which is the region, shown in Fig. \ref{RocF2} (\textit{Left}),
where the $F_2$ series coincides with the Appell $F_2$ function. Outside of this region, the Appell $F_2$ function can be defined by the integral representation of the Euler type
\begin{align}
&F_2(a,b_1,b_2;c_1,c_2;x,y)= 
\nonumber\\
&\frac{\Gamma(c_1)\Gamma(c_2)}{\Gamma(b_1)\Gamma(b_2)\Gamma(c_1-b_1)\Gamma(c_2-b_2)}\int_0^1du\int_0^1\ dv\ u^{b_1-1}v^{b_2-1}(1-u)^{c_1-b_1-1}(1-v)^{c_2-b_2-1}(1-ux-vy)^{-a}\label{F_2_Euler}
\end{align}
subject to the constraints that $\text{Re}(b_1),\text{Re}(b_2),\text{Re}(c_1-b_1)$ and $\text{Re}(c_2-b_2)$ are positive numbers, or by 
\begin{align}
 F_2(a,b_1,b_2;c_1,c_2;x,y)= \frac{\Gamma (c_2)}{\Gamma (b_2) \Gamma (c_2-b_2)}  \int_0^1 d v v^{b_2-1}  (1-v y)^{-a} (1-v)^{-b_2+c_2-1}  {}_2F_1\left(a,b_1;c_1;\frac{x}{1-v y}\right)\label{F_2_Euler2}
\end{align}

Another well-known integral representation of $F_2$ is of the Mellin-Barnes type
\begin{align}
F_2(a,b_1,b_2;&c_1,c_2;x,y)= 
\nonumber\\
&\frac{\Gamma(c_1)\Gamma(c_2)}{\Gamma(a)\Gamma(b_1)\Gamma(b_2)}\int_{-i\infty}^{+i\infty} ds\int_{-i\infty}^{+i\infty} dt\ (-x)^s(-y)^t\Gamma(-s)\Gamma(-t)\frac{\Gamma(a+s+t)\Gamma(b_1+s)\Gamma(b_2+t)}{\Gamma(c_1+s)\Gamma(c_2+t)}\label{F_2_MB}
\end{align}
where the integration contours are such that they separate the poles of $\Gamma(-s)$ and $\Gamma(-t)$ from those of $\Gamma(a+s+t)$, $\Gamma(b_1+s)$ and $\Gamma(b_2+t)$.

 $F_2$ has the following symmetry
\begin{equation}\label{F_2sym}
F_2(a,b_1,b_2;c_1,c_2;x,y)= F_2(a,b_2,b_1;c_2,c_1;y,x)
\end{equation}
and from suitable changes of variables in Eq.(\ref{F_2_Euler}), one can obtain its well-known Euler transformations \cite{Appell}
\begin{align}
F_2(a,b_1,b_2;c_1,c_2;x,y)&=(1-x)^{-a}F_2\left(a,c_1-b_1,b_2;c_1,c_2;\frac{x}{x-1},\frac{y}{1-x}\right)\nonumber\\
&=(1-y)^{-a}F_2\left(a,b_1,c_2-b_2;c_1,c_2;\frac{x}{1-y},\frac{y}{y-1}\right)\nonumber\\
&=(1-x-y)^{-a}F_2\left(a,c_1-b_1,c_2-b_2;c_1,c_2;\frac{x}{x+y-1},\frac{y}{x+y-1}\right)
 \label{F_2_EulerTransforms}
\end{align}
which will be useful in the following.

The system of partial differential equations satisfied by $F_2$ is given by \cite{Appell}
\begin{align}
x(1-x)r - xy s + [c_{1}-(a+b_{1}+1)x]p-b_{1} y q-a b_{1} z =0 \nonumber\\
y(1-y)t - xy s + [c_{2}-(a+b_{2}+1)y]q-b_{2} x p-a b_{2} z =0\label{PDE}
\end{align}
where $r=z_{xx}$, $t=z_{yy}$, $s= z_{xy}$, $p=z_{x}$, $q=z_{y}$.

The singular curves of the above system are  
$x=0$, $y=0$, $x=1$, $y=1$, $x+y=1$. They are shown in Fig. \ref{RocF2} (\textit{Right}).

We will now consider analytic continuations of the Appell $F_2$ series with the aim to evaluate it for generic values of its parameters $a, b_1, b_2, c_1, c_2$ and for all possible real values of $x$ and $y$ except on these singular curves. 

\section{A first analytic continuation study based on the Mellin-Barnes representation of $F_2$
\label{MB_F2Sec}}

It is straightforward to derive two well-known analytic continuation formulas (and two related symmetrical expressions) of the Appell $F_2$ series from the Mellin-Barnes representation presented in Eq.(\ref{F_2_MB}). For this, one can use the method of \cite{ABFGgeneral} or, equivalently, of \cite{Friot:2011ic, Passare:1996db, TZ}, which give 
 \begin{align}\label{F2_ACMB1}
 F_2(a,b_1,b_2;c_1,c_2;x,y)= \frac{\g{c_2}\g{b_2-a}}{\g{b_2}\g{c_2-a}}(-y)^{-a}F{}^{2:1;0}_{1:1;0}
  \left[
   \setlength{\arraycolsep}{0pt}
   \begin{array}{c@{{}:{}}c@{;{}}c}
 a,a-c_2+1 & b_1&-\\[1ex]
  a-b_2+1 & c_1 & \linefill
   \end{array}
   \;\middle|\;
 -\frac{x}{y},\frac{1}{y}
 \right]\\\nonumber
 + \frac{\g{c_2}\g{a-b_2}}{\g{a}\g{c_2-b_2}}(-y)^{-b_2}\,H_2\left(a-b_2,b_1,b_2,b_2-c_2+1;c_1;x,-\frac{1}{y}\right)
 \end{align}
and
\begin{multline}\label{F2_ACMB2}
 F_2(a,b_1,b_2,c_1,c_2;x,y)= \frac{\g{c_2}\g{b_2-a}}{\g{b_2}\g{c_2-a}}(-y)^{-a} F{}^{2:1;0}_{1:1;0}
  \left[
   \setlength{\arraycolsep}{0pt}
   \begin{array}{c@{{}:{}}c@{;{}}c}
 a,a-c_2+1 & b_1& \linefill\\[1ex]
  a-b_2+1 & c_1 & \linefill
   \end{array}
   \;\middle|\;
-\frac{x}{y},\frac{1}{y}
 \right]\\
 +\frac{\g{c_2}\g{a-b_2}\g{c_1}\g{b_1+b_2-a}}{\g{a}\g{c_2-b_2}\g{b_1}\g{c_1+b_2-a}}(-x)^{b_2-a}(-y)^{-b_2}
 \tilde{F}{}^{2:0;2}_{1:0;0}
  \left[
   \setlength{\arraycolsep}{0pt}
   \begin{array}{c@{{}:{}}c@{;{}}c}
 a-b_2 ,a-b_2-c_1+1& \linefill &b_2,b_2-c_2+1\\[1ex]
  a-b_1-b_2+1& - & -
   \end{array}
   \;\middle|\;
\frac{1}{x},\frac{x}{y}
 \right]\\
 + \frac{\g{c_1}\g{c_2}\g{a-b_1-b_2}}{\g{a}\g{c_2-b_2}\g{c_1-b_1}}(-x)^{-b_1}(-y)^{-b_2}
 F_{3}\Big(b_1,b_2,b_1-c_1+1,b_2-c_2+1,b_1+b_2-a+1;\frac{1}{x},\frac{1}{y}\Big)
 \end{multline}
 where a commonly used notation for the Kamp\'e de F\'eriet series is, with $(a_p)\doteq a_1,...,a_p$, \cite{Srivastava}
\begin{align}\label{KdFnotation}
F{}^{p:q;k}_{l:m;n}
  \left[
   \setlength{\arraycolsep}{0pt}
   \begin{array}{c@{{}:{}}c@{;{}}c}
  (a_p) & (b_q) & (c_k)\\[1ex]
   (\alpha_l) & (\beta_m) & (\gamma_n)
   \end{array}
   \;\middle|\;
 x,y
 \right]=
 \sum_{r=0}^{\infty}\sum_{s=0}^{\infty}\frac{\prod_{j_1=1}^{p}(a_{j_1})_{r+s}\prod_{j_2=1}^{q}(b_{j_2})_{r}\prod_{j_3=1}^{k}(c_{j_3})_{s}}{\prod_{j_4=1}^{l}(\alpha_{j_4})_{r+s}\prod_{j_5=1}^{m}(\beta_{j_5})_{r}\prod_{j_6=1}^{n}(\gamma_{j_6})_{s}}\frac{x^r}{r!}\frac{y^s}{s!}
\end{align} 
and where the $\tilde{F}$ double series is defined as
\begin{align}\label{MirrorKdFnotation}
\tilde{F}{}^{p:q;k}_{l:m;n}
  \left[
   \setlength{\arraycolsep}{0pt}
   \begin{array}{c@{{}:{}}c@{;{}}c}
  (a_p) & (b_q) & (c_k)\\[1ex]
   (\alpha_l) & (\beta_m) & (\gamma_n)
   \end{array}
   \;\middle|\;
 x,y
 \right]\doteq
 \sum_{r=0}^{\infty}\sum_{s=0}^{\infty}\frac{\prod_{j_1=1}^{p}(a_{j_1})_{r-s}\prod_{j_2=1}^{q}(b_{j_2})_{r}\prod_{j_3=1}^{k}(c_{j_3})_{s}}{\prod_{j_4=1}^{l}(\alpha_{j_4})_{r-s}\prod_{j_5=1}^{m}(\beta_{j_5})_{r}\prod_{j_6=1}^{n}(\gamma_{j_6})_{s}}\frac{x^r}{r!}\frac{y^s}{s!}
\end{align} 
One will note that, as performed in Appendix C-2 of \cite{DelDuca:2009ac}, $\tilde{F}$ can be transformed in terms of Kamp\'e de F\'eriet series if necessary.

Eqs.(\ref{F2_ACMB1}) and (\ref{F2_ACMB2}) respectively match with Eq.(64) p.294 and Eq.(66) p.295 of \cite{Srivastava}.

It is easy to derive the convergence regions of these analytic continuations from the well-known convergence properties of the Appell $F_3$, Horn $H_2$ and Kamp\'e de F\'eriet series, and by noting that, from the property of cancellation of opposite elements in the characteristic list of a hypergeometric series\footnote{See Section \ref{FirstAC01} for a brief reminder about this fact.}, $\tilde{F}$ in Eq.(\ref{F2_ACMB2}) has the same convergence properties as $H_2$ with the same arguments.

\begin{figure}[h]
\centering
\includegraphics[width=6cm, height=6cm]{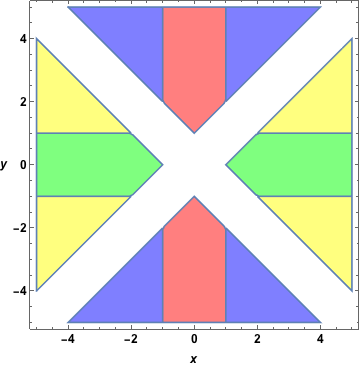}
\includegraphics[width=6cm, height=6cm]{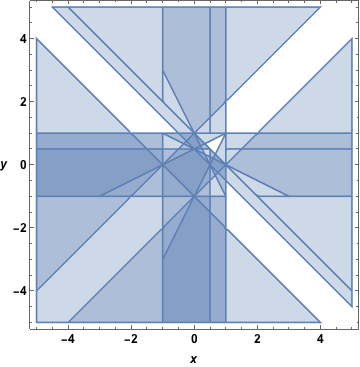}
\caption{\textit{Left:} Regions of convergence of the analytic continuations of the Appell $F_2$ series given in Eqs.(\ref{F2_ACMB1}) (in red) and (\ref{F2_ACMB1}) (in blue) and of the symmetrical relations obtained from them (in green and yellow respectively). \textit{Right:} Same figure when one adds the ROC of Fig. \ref{RocF2} (\textit{Left}) and when one applies the Euler transforms of $F_2$ given in Eq. (\ref{F_2_EulerTransforms}) on all these results.\label{RocF2MB}}
\end{figure}

One then obtains for Eq.(\ref{F2_ACMB1}) the convergence region $ |x|<1 \wedge |-\frac{1}{y}|<1\wedge |-\frac{1}{y}|(1+|x|)<1 $ shown in red in Fig. \ref{RocF2MB} (\textit{Left}), for real values of $x$ and $y$. As for Eq.(\ref{F2_ACMB2}), the convergence region is $|\frac{1}{x}|<1 \wedge |\frac{x}{y}|<1 \wedge |\frac{x}{y}|(1+ |\frac{1}{x}|)<1 $,
and it is shown in blue in the same figure. 

As mentioned above, two symmetrical relations can be computed from the Mellin-Barnes representation, which can also be obtained using the symmetry property of Eq.(\ref{F_2sym}) applied to Eq.(\ref{F2_ACMB1}) and Eq.(\ref{F2_ACMB2}). These symmetrical analytic continuations converge in the green and yellow regions of Fig. \ref{RocF2MB} (\textit{Left}).

With no further transformation of the MB integral one cannot obtain, from the latter, other series representations than those presented above. However, using the three Euler transformations shown in Eq.(\ref{F_2_EulerTransforms}), it is possible to derive 12 other formulas, that we do not list here and which, alltogether, allow us to obtain the total convergence region of Fig. \ref{RocF2MB} (\textit{Right}).

When added to the usual $F_2$ series definition and its three Euler transformations, these 16 linear transformations show that a good part of the $(x,y)$ real plane can be reached, but one can see on the plot that several regions are still missing. This is the aim of the rest of this paper to fill this gap, following an alternative method of analytic continuation.

\section{Analytic continuation from Olsson's method\label{Olsson_Sec}}

In \cite{Olsson-64}, Olsson obtained the complete set of solutions of the Appell $F_1$ system of partial differential equations, as well as the relations that connect them, thereby obtaining analytic continuation formulas of $F_1$ in the neighbourhood of any of its singularities. His method rests on the application of various transformations and analytic continuations of the $_2F_1$ Gauss hypergeometric series on the Appell $F_1$ written as an infinite sum of $_2F_1$. We follow this procedure below to derive analytic continuation formulas for $F_2$.

One will note that, except for the results of Section \ref{MB_F2Sec}, which are briefly rederived following Olsson's method in the beginning of subsection \ref{OlssonMB}, the regions of convergence of all the analytic continuation formulas presented in Section \ref{Olsson_Sec} are trivial. They can be straightforwardly obtained as the intersections of the regions defined by the modulus, smaller than unity, of each of the arguments of the series involved in these formulas. This is due to the simple form of these series, as it is explictly shown in one example in Section \ref{FirstAC01}.

\subsection{Analytic continuation around (0,1)\label{0,1} }
We begin our study by the derivation of analytic continuations of the Appell $F_2$ series around the point $(0,1)$. Note that, still by the symmetry shown in Eq. (\ref{F_2sym}), the final expressions can be used to obtain analytic continuations around the point $(1,0)$. Several different formulas will be necessary to cover the whole neighbourhood of these points.

 \subsubsection{A first analytic continuation\label{FirstAC01}}
Rewriting $F_2$ as an infinite sum of $_2F_1$, one gets
\begin{equation}\label{F2as2F1}
F_2(a,b_1,b_2;c_1,c_2;x,y)= \sum_{m=0}^{\infty}\frac{  \p{a}{m}\p{b_1}{m}}{ \p{c_1}{m} m!}x^m \,_2F_1(a+m,b_2;c_2;y)
\end{equation}
where one can now use the well-known analytic continuation of \(_2F_1(a,b,c;z)\) around \(z=1\) given by 
\begin{align}
 _2 F_1(a,b,c;z)&= \frac{\Gamma(c)\Gamma(c-a-b)}{\Gamma(c-a)\Gamma(c-b)}\,_2 F_1(a,b,a+b-c+1;1-z)\nonumber\\
 &+ \frac{\Gamma(c)\Gamma(a+b-c)}{\Gamma(a)\Gamma(b)}(1-z)^{c-a-b}\,_2 F_1(c-a,c-b,c-a-b+1;1-z)\label{ana_1}
 \end{align}
Substituting and simplifying, one obtains
 \begin{align}
 F_2(a,b_1,b_2;c_1,c_2;x,y)= \frac{\g{c_2}\g{c_2 - b_2-a}}{\g{c_2-a}\g{c_2-b_2}}F{}^{1:2;1}_{1:1;0}
  \left[
   \setlength{\arraycolsep}{0pt}
   \begin{array}{c@{{}:{}}c@{;{}}c}
 a & b_1,1+a-c_2& b_2\\[1ex]
  a+b_2-c_2+1 & c_1 & \linefill
   \end{array}
   \;\middle|\;
 x,1-y
 \right]\nonumber\\
 + \frac{\g{c_2}\g{a+b_2-c_2 }}{\g{a}\g{b_2}}(1-y)^{c_2-b_2-a}\,\tilde{F}{}^{1:1;2}_{1:0;1}
  \left[
   \setlength{\arraycolsep}{0pt}
   \begin{array}{c@{{}:{}}c@{;{}}c}
 c_2-a & c_2-b_2 & b_1,1+a-c_2\\[1ex]
  1+c_2-b_2-a & \linefill & c_1
   \end{array}
   \;\middle|\;
1-y,\frac{x}{1-y}
 \right]
 \label{F_2sol011}
 \end{align}
 where we have used, for the Kamp\'e de F\'eriet and $\tilde{F}$ series, the notation presented in the previous section.
 
As mentioned above, a simple look at the particular form of these series allows to straightforwardly conclude that the $F^{1:2;1}_{1:1;0}$ series converges for $ |x|<1 \wedge |1-y|<1 $ (see Fig. \ref{Fig(0,1)} \textit{Left}) 
and the $\tilde{F}{}^{1:1;2}_{1:0;1}$ series converges for $|\frac{x}{1-y}|<1 \wedge |1-y|<1 $ (see Fig. \ref{Fig(0,1)} \textit{Right}). This is due to the fact that the convergence region of a hypergeometric series is independent of its parameters (exceptional values of the latter being excluded) which implies, in particular, that the cancellation of opposite elements in the characteristic list of this series does not affect its region of convergence \cite{Srivastava}. In our present case of study, since the characteristic lists of the Kamp\'e de F\'eriet and $\tilde{F}{}^{1:1;2}_{1:0;1}$ series are respectively $\{m+n,m,m,n,-(m+n),-m\}$ and $\{m-n,m,n,-(m-n),-n\}$, the cancellation property leads to a factorisation of these double series into single series whose convergence regions are trivial.  
\begin{figure}[h]
\begin{center}
  \includegraphics[scale=.24]{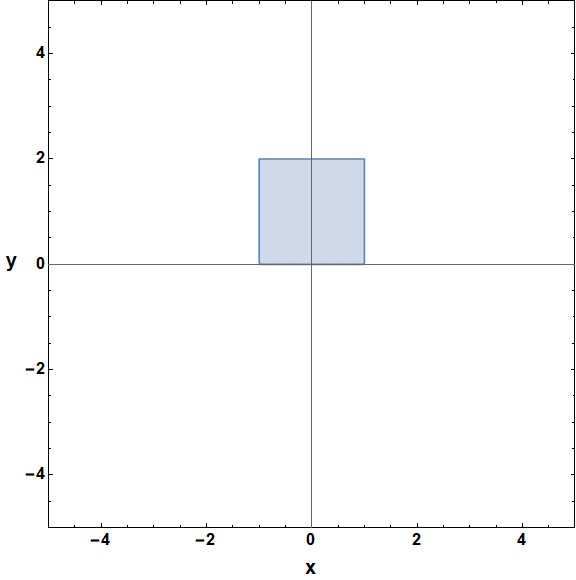}
  \includegraphics[scale=.24]{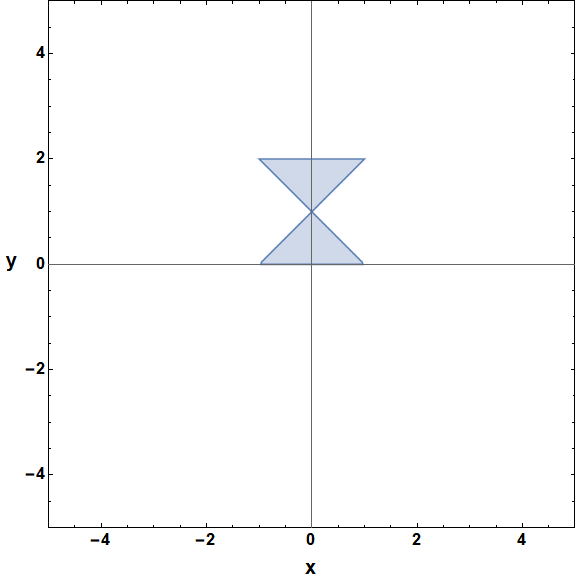}
\caption{\textit{Left}: ROC of $F^{1:2;1}_{1:1;0}$ in Eq.\eqref{F_2sol011}. \textit{Right}: ROC of $\tilde{F}$ in Eq.\eqref{F_2sol011}.\label{Fig(0,1)}}
\end{center}
\end{figure}\\\\
Therefore, from Fig.\ref{Fig(0,1)} one concludes that, for real values of $x$ and $y$, the ROC of the RHS of Eq. \eqref{F_2sol011} is restricted to the region shown in Fig. \ref{Fig(0,1)} \textit{Right}.\\
Using Eq.\eqref{F_2sym} we get the analytic continuation around (1,0) as
\begin{multline}
 F_2(a,b_1,b_2;c_1,c_2;x,y)= \frac{\g{c_1}\g{c_1 - b_1-a}}{\g{c_1-a}\g{c_1-b_1}}F{}^{1:2;1}_{1:1;0}
  \left[
   \setlength{\arraycolsep}{0pt}
   \begin{array}{c@{{}:{}}c@{;{}}c}
 a & b_2,1+a-c_1& b_1\\[1ex]
  a+b_1-c_1+1 & c_2 & \linefill
   \end{array}
   \;\middle|\;
 y,1-x
 \right]\\\\
 + \frac{\g{c_1}\g{a+b_1-c_1 }}{\g{a}\g{b_1}}(1-x)^{c_1-b_1-a}\tilde{F}{}^{1:1;2}_{1:0;1}
  \left[
   \setlength{\arraycolsep}{0pt}
   \begin{array}{c@{{}:{}}c@{;{}}c}
 c_1-a & c_1-b_1 & b_2,1+a-c_1\\[1ex]
  1+c_1-b_1-a & \linefill & c_2
   \end{array}
   \;\middle|\;
1-x,\frac{y}{1-x}
 \right]\label{F_2sol101}
 \end{multline}
which now converges, still for real values of $x$ and $y$, in the region shown in Fig.\ref{F_210af_t}. One will note here that the regions shown in Fig.\ref{Fig(0,1)} \textit{Right} and Fig.\ref{F_210af_t} are already covered by the formulas presented in Section \ref{MB_F2Sec}.
\begin{figure}[h]
  \centering
  \includegraphics[scale=.24]{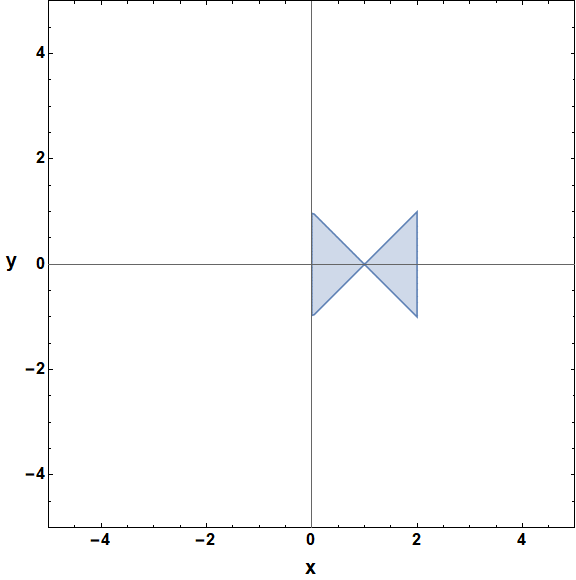}
  \caption{ROC of the RHS of Eq.\eqref{F_2sol101}.}\label{F_210af_t}
\end{figure}
Looking at the ROC plots it is clear that the $\tilde{F}$ series prevents the analytic continuation in Eq.\eqref{F_2sol011} (respectively Eq.\eqref{F_2sol101}) to converge around the whole neighbourhood of $(0,1)$ (respectively $(1,0)$). More precisely, the constraint $|\frac{x}{1-y}|<1$ (respectively $|\frac{y}{1-x}|<1$) is responsible for that, which means that by going on with the analytic continuation process on its associated sum, one can probably derive another more interesting formula.
This will be done in the next section where the obtained formula, in addition to reach the missing region around (0,1), will also be the starting point of another analytic continuation which is presented in Section \ref{inf,1}.

 \subsubsection{A second analytic continuation}

\begin{figure}[h]
\centering
  \includegraphics[scale=.24]{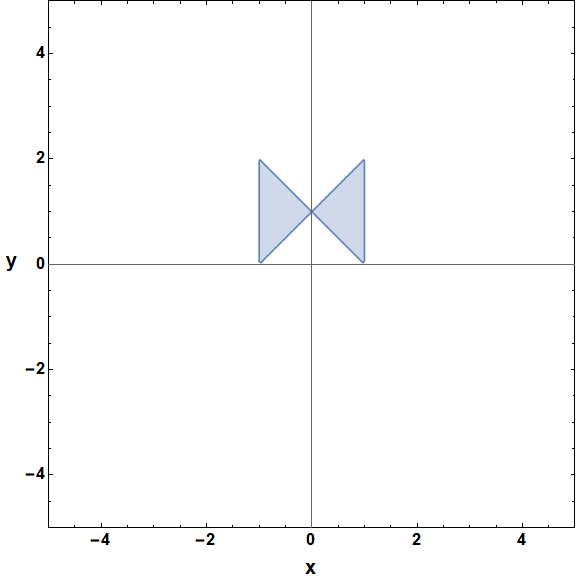}
   \includegraphics[scale=.24]{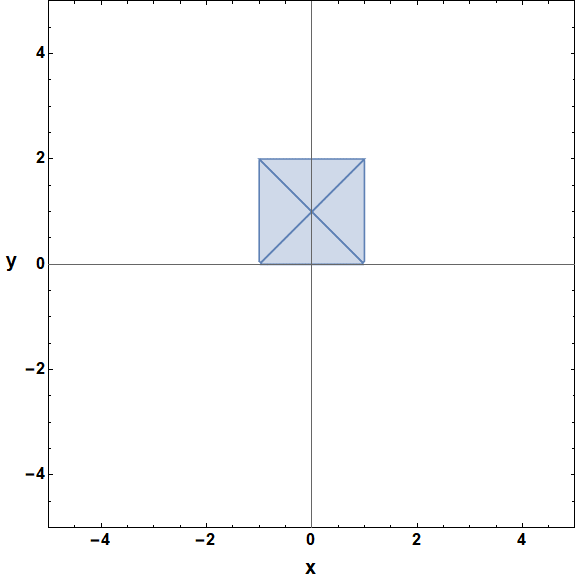}
   \includegraphics[scale=.24]{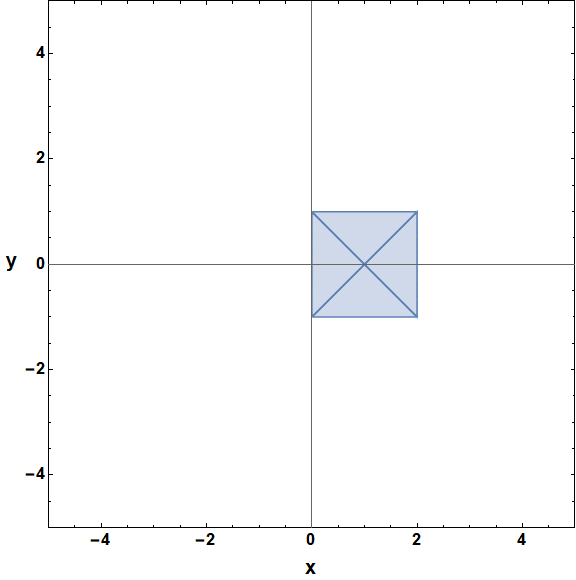}
  \caption{\textit{Left:} ROC of the RHS of Eq.\eqref{F_2sol012}, \textit{Middle:} Reached region around $(0,1)$, \textit{Right:} Corresponding reached region around $(1,0)$. \label{F_201bfinal}}
\end{figure}

Since we have
\begin{multline}
\tilde{F}{}^{1:1;2}_{1:0;1}
  \left[
   \setlength{\arraycolsep}{0pt}
   \begin{array}{c@{{}:{}}c@{;{}}c}
 c_2-a & c_2-b_2 & b_1,1+a-c_2\\[1ex]
  1+c_2-b_2-a & \linefill & c_1
   \end{array}
   \;\middle|\;
1-y,\frac{x}{1-y}
 \right]= \\
 \sum_{m=0}^{\infty}\sum_{n=0}^{\infty}\frac{\p{c_2 -a}{m-n}\p{b_1}{n}\p{1+a-c_2}{n}\p{c_2-b_2}{m}}{\p{1+c_2-b_2-a}{m-n}\p{c_1}{n}m!n!}(1-y)^{m}\Big(\frac{x}{1-y}\Big)^{n}
\end{multline}
which can be rewritten as 
\begin{multline}
\tilde{F}{}^{1:1;2}_{1:0;1}
  \left[
   \setlength{\arraycolsep}{0pt}
   \begin{array}{c@{{}:{}}c@{;{}}c}
 c_2-a & c_2-b_2 & b_1,1+a-c_2\\[1ex]
  1+c_2-b_2-a & \linefill & c_1
   \end{array}
   \;\middle|\;
1-y,\frac{x}{1-y}
 \right]= \\
 \sum_{m=0}^{\infty} \frac{\p{c_2 -a}{m}\p{c_2-b_2}{m}}{\p{1+c_2-b_2-a}{m}m!}(1-y)^{m}\,_3F_2
\left(
\begin{matrix}
a+b_2-c_2-m, & b_1, & 1+a-c_2\\\\
&\hspace{-.8cm}1+a-c_2-m, &  c_1   
\end{matrix}
\Bigg {|} { \frac{x}{1-y} }
\right)
\end{multline}
one can use the standard analytic continuation formula of $_3F_2$ 
\begin{equation}\label{3F2ac}
\begin{aligned}
_3F_2
\left(
\begin{matrix}
a_1, & a_2, & a_3 \\\\
&\hspace{-.3cm}b_1, & \hspace{-.3cm}b_2   
\end{matrix}
\Bigg {|} { z }
\right)&= \frac{\g{b_1}\g{b_2}\g{a_2-a_1}\g{a_3-a_1}}{\g{a_2}\g{a_3}\g{b_1-a_1}\g{b_2-a_1}}(-z)^{-a_1} \,_3F_2
\left(
\begin{matrix}
a_1, &1+a_1-b_1, & 1+a_1-b_2 \\\\
&\hspace{-.3cm}1+a_1-a_2, & \hspace{-.3cm}1+a_1-a_3   
\end{matrix}
\Bigg {|} { \frac{1}{z}}
\right)\\\\
&+ \frac{\g{b_1}\g{b_2}\g{a_1-a_2}\g{a_3-a_2}}{\g{a_1}\g{a_3}\g{b_1-a_2}\g{b_2-a_2}}(-z)^{-a_2}\,_3F_2
\left(
\begin{matrix}
a_2, &1+a_2-b_1, & 1+a_2-b_2 \\\\
&\hspace{-.3cm}1+a_2-a_1, & \hspace{-.3cm}1+a_2-a_3   
\end{matrix}
\Bigg {|} { \frac{1}{z}}
\right) \\\\
&+\frac{\g{b_1}\g{b_2}\g{a_1-a_3}\g{a_2-a_3}}{\g{a_1}\g{a_2}\g{b_1-a_3}\g{b_2-a_3}}(-z)^{-a_3}\,_3F_2
\left(
\begin{matrix}
a_3, & 1+a_3-b_1, & 1+a_3-b_2 \\\\
&\hspace{-.3cm}1+a_3-a_1, & \hspace{-.3cm}1+a_3-a_2   
\end{matrix}
\Bigg {|} {\frac{1}{z}}
\right)
\end{aligned}
\end{equation}
which allows to derive 
\begin{align}
    &\tilde{F}{}^{1:1;2}_{1:0;1}
  \left[
   \setlength{\arraycolsep}{0pt}
   \begin{array}{c@{{}:{}}c@{;{}}c}
 c_2-a & c_2-b_2 & b_1,1+a-c_2\\[1ex]
  1+c_2-b_2-a & \linefill & c_1
   \end{array}
   \;\middle|\;
1-y,\frac{x}{1-y}
 \right]= \nonumber\\ &\left(\frac{x}{y-1}\right)^{-a-b_2+c_2} \frac{\Gamma \left(c_1\right) \Gamma \left(-a+b_1-b_2+c_2\right) }{\Gamma \left(b_1\right) \Gamma \left(-a-b_2+c_1+c_2\right)} \tilde{F}{}^{2:1;1}_{2:0;0}
  \left[
   \setlength{\arraycolsep}{0pt}
   \begin{array}{c@{{}:{}}c@{;{}}c}
 1-b_2 , -a+b_1-b_2+c_2 & c_2-b_2 & b_2  \\[1ex]
  -a-b_2+c_2+1,-a-b_2+c_1+c_2 &\linefill & \linefill
   \end{array}
   \;\middle|\;
 x,\frac{1-y}{x}
 \right]\nonumber\\
 &+\left(\frac{x}{y-1}\right)^{-b_1} \frac{\Gamma \left(c_1\right)  \Gamma \left(a-b_1+b_2-c_2\right)}{\Gamma \left(c_1-b_1\right) \Gamma \left(a+b_2-c_2\right)}F{}^{1:1;2}_{1:0;1}
  \left[
   \setlength{\arraycolsep}{0pt}
   \begin{array}{c@{{}:{}}c@{;{}}c}
 -a+b_1+c_2 & c_2-b_2& b_1-c_1+1, b_1\\[1ex]
  -a+b_1-b_2+c_2+1 & \linefill &-a+b_1+c_2
   \end{array}
   \;\middle|\;
 1-y,\frac{1-y}{x}
 \right]
\end{align}
Substituting this result back in Eq.\eqref{F_2sol011} one gets 
\begin{align}\label{F_2sol012}
F_2(a,b_1,b_2;c_1,c_2;x,y)
= &\frac{\Gamma (c_2) \Gamma (-a-b_2+c_2)}{\Gamma (c_2-a) \Gamma (c_2-b_2)} {F}{}^{1:2;1}_{1:1;0}
  \left[
   \setlength{\arraycolsep}{0pt}
   \begin{array}{c@{{}:{}}c@{;{}}c}
  a & a-c_2+1,b_1 & b_2\\[1ex]
   a+b_2-c_2+1 & c_1 & -
   \end{array}
   \;\middle|\;
 x,1-y
 \right]\nonumber\\
		&+(1-y)^{-a-b_2+c_2}\left(\frac{x}{y-1}\right)^{-a-b_2+c_2} \frac{\Gamma (c_1) \Gamma (c_2)  \Gamma (a+b_2-c_2) \Gamma (-a+b_1-b_2+c_2) }{\Gamma (a) \Gamma (b_1) \Gamma (b_2) \Gamma (-a-b_2+c_1+c_2)}\nonumber\\
		&\times\tilde{F}{}^{2:1;1}_{2:0;0}
  \left[
   \setlength{\arraycolsep}{0pt}
   \begin{array}{c@{{}:{}}c@{;{}}c}
  a+b_2-c_2, a+b_2-c_2-c_1+1 & b_2 & c_2-b_2\\[1ex]
   b_2,a-b_1+b_2-c_2+1 & - & -
   \end{array}
   \;\middle|\;
 \frac{1-y}{x},x
 \right]\nonumber\\
		&+\left(\frac{x}{y-1}\right)^{-b_1} (1-y)^{-a-b_2+c_2}\frac{\Gamma (c_1) \Gamma (c_2)  \Gamma (a-b_1+b_2-c_2)}{\Gamma (a) \Gamma (b_2) \Gamma (c_1-b_1)}\nonumber\\
		&\times{F}{}^{1:2;1}_{1:1;0}
  \left[
   \setlength{\arraycolsep}{0pt}
   \begin{array}{c@{{}:{}}c@{;{}}c}
  -a+b_1+c_2 & b_1,b_1-c_1+1 & c_2-b_2\\[1ex]
   -a+b_1-b_2+c_2+1 & -a+b_1+c_2 & -
   \end{array}
   \;\middle|\;
 \frac{1-y}{x},1-y
 \right]
	\end{align}
which, trivially, converges in $ |\frac{1-y}{x}|<1 \wedge |x|<1  \wedge |1-y|<1$ (see, for real values of $x$ and $y$, Fig. \ref{F_201bfinal} \textit{Left}).

Therefore, taking into account the results of Section \ref{FirstAC01}, we are now able to reach the whole neighborhood of $(0,1)$ (see Fig.\ref{F_201bfinal} \textit{Middle}) with the exception of the boundaries of the ROCs which will be considered later.

As before we can now straightforwardly get the corresponding continuation around $(1,0)$ using Eq.\eqref{F_2sym}. Together with the previous symmetrical analytic continuation given in Eq.(\ref{F_2sol101}) we can then reach the whole neighborhood of $(1,0)$, except on the boundaries of the ROCs (see Fig.\ref{F_201bfinal} \textit{Right}).

\subsection{Analytic continuation around $(0,\infty)$ and $(\infty,\infty)$\label{OlssonMB}}
Next we turn to the analytic continuation around the singular point $(0,\infty)$. This analytic continuation will then be used to find the analytic continuation around the singular point $(\infty,\infty)$. Moreover, as before, symmetry will give us the $(\infty,0)$ case.
The corresponding formulas and their convergence regions have already been derived in Section \ref{MB_F2Sec}, from the Mellin-Barnes representation of $F_2$, but for completeness we briefly show here how they can also be obtained from Olsson's method (see also Chapter 9 of \cite{Srivastava}). 

To find the continuation around $(0,\infty)$ we start again with Eq.\eqref{F2as2F1} and use the standard analytic continuation formula for $_2F_1$  given by
\begin{align}
 _2 F_1(a,b,c;z)&= \frac{\Gamma(c)\Gamma(b-a)}{\Gamma(b)\Gamma(c-a)}(-z)^{-a}\,_2 F_1\left(a,a-c+1,a-b+1;\frac{1}{z}\right)\nonumber\\
 &+ \frac{\Gamma(c)\Gamma(a-b)}{\Gamma(a)\Gamma(c-b)}(-z)^{-b}\,_2 F_1\left(b,b-c+1,b-a+1;\frac{1}{z}\right)\label{ana_2}
 \end{align}
This directly gives Eq.(\ref{F2_ACMB1}).

Then, in order to find the analytic continuation around $(\infty,\infty)$, we observe that we just need to continue the  Horn $H_2$ series in Eq.(\ref{F2_ACMB1}), as the other series already converges around $(\infty,\infty)$.
So, writing 
 \begin{equation}
H_{2}\Big(a-b_2,b_1,b_2,b_2-c_2+1;c_1;x,-\frac{1}{y}\Big)=\sum_{n=0}^{\infty} \frac{\p{a-b_2}{-n}\p{b_2}{n}\p{b_2-c_2+1}{n}}{n!}\Big(-\frac{1}{y}\Big)^{n}\,_2F_1(a-b_2-n,b_1,c_1;x)
 \end{equation}
and using once more Eq.\eqref{ana_2} we get
 \begin{multline}
H_{2}\Big(a-b_2,b_1,b_2,b_2-c_2+1;c_1;x,-\frac{1}{y}\Big)=\\
\frac{\g{c_1}\g{b_1+b_2-a}}{\g{b_1}\g{c_1+b_2-a}}(-x)^{b_2-a}
\tilde{F}{}^{2:0;2}_{1:0;0}
  \left[
   \setlength{\arraycolsep}{0pt}
   \begin{array}{c@{{}:{}}c@{;{}}c}
 a-b_2 ,a-b_2-c_1+1& \linefill &b_2,b_2-c_2+1\\[1ex]
  a-b_2-b_1+1&\linefill & \linefill
   \end{array}
   \;\middle|\;
 \frac{1}{x},\frac{x}{y}
 \right]\\\\
 + \frac{\g{c_1}\g{a-b_1-b_2}}{\g{a-b_2}\g{c_1-b_1}}(-x)^{-b_1}
  F_{3}\Big(b_1,b_2,b_1-c_1+1,b_2-c_2+1,b_1+b_2-a+1;\frac{1}{x},\frac{1}{y}\Big)
 \end{multline}
Substituting the above back in Eq.(\ref{F2_ACMB1}) we get Eq.(\ref{F2_ACMB2}).

The regions of convergence of Eqs.(\ref{F2_ACMB1}) and (\ref{F2_ACMB2}) are ploted in Fig. \ref{RocF2MB} \textit{Left}.

\subsection{Analytic Continuation around $(\infty,1)$\label{inf,1}}
Looking at the ROC of the second continuation around (0,1), Eq.\eqref{F_2sol012}, we note that $F{}^{1:2;1}_{1:1;0}(...;\frac{1-y}{x},1-y)$ converges in the entire neighborhood of $(\infty,1)$. 
Therefore, in order to derive an analytic continuation valid in this region, we only need to analytically continue the remaining two series of Eq.\eqref{F_2sol012}.
Let us consider the first series which reads
\begin{multline}
F{}^{1:2;1}_{1:1;0}
  \left[
   \setlength{\arraycolsep}{0pt}
   \begin{array}{c@{{}:{}}c@{;{}}c}
 a & b_1,1+a-c_2& b_2\\[1ex]
  a+b_2-c_2+1 & c_1 & \linefill
   \end{array}
   \;\middle|\;
 x,1-y
 \right]=\\\\
  \sum_{n=0}^{\infty}\frac{\p{a}{n}\p{b_2}{n}}{\p{a+b_2-c_2+1}{n}n!}(1-y)^{n}
\, _3F_2
\left(
\begin{matrix}
a+n, &\hspace{-.7cm} b_1, &\hspace{-.7cm} 1+a-c_2\\\\
&\hspace{-.3cm}a+b_2-c_2+1+n, & \hspace{-.3cm}c_1   
\end{matrix}
\Bigg {|} {x}
\right) 
\end{multline}
Applying Eq.\eqref{3F2ac} and simplifying we get
\begin{align}\label{F2_1infa}
&F{}^{1:2;1}_{1:1;0}
  \left[
   \setlength{\arraycolsep}{0pt}
   \begin{array}{c@{{}:{}}c@{;{}}c}
 a & b_1,1+a-c_2& b_2\\[1ex]
  a+b_2-c_2+1 & c_1 & \linefill
   \end{array}
   \;\middle|\;
 x,1-y
 \right]= \nonumber\\
 &\frac{\g{1+a-c_2-b_1}\g{c_1}\g{a-b_1}\g{a+b_2-c_2+1}}{\g{1+a-c_2}\g{c_1-b_1}\g{a}\g{a+b_2-b_1-c_2+1}}(-x)^{b_1}
 \tilde{F}{}^{1:2;1}_{1:1;0}
  \left[
   \setlength{\arraycolsep}{0pt}
   \begin{array}{c@{{}:{}}c@{;{}}c}
 b_1+c_2-a-b_2 & 1+b_1-c_1,b_1& b_2\\[1ex]
  1-a+b_1 & b_1+c_2-a & \linefill
   \end{array}
   \;\middle|\;
 \frac{1}{x},1-y
 \right]\nonumber\\
 &+ \frac{\g{b_1+c_2-a-1}\g{c_1}\g{c_2-1}\g{a+b_2-c_2+1}}{\g{b_1}\g{c_1+c_2-a-1}\g{a}\g{b_2}}(-x)^{c_2-a-1}\nonumber\\
 &\hspace{3cm}\times\tilde{F}{}^{1:2;1}_{1:1;0}
  \left[
   \setlength{\arraycolsep}{0pt}
   \begin{array}{c@{{}:{}}c@{;{}}c}
 1-b_2 & 2+a-c_2-c_1,1+a-c_2& b_2\\[1ex]
  2-c_2 & a-c_2-b_1+2 & \linefill
   \end{array}
   \;\middle|\;
 \frac{1}{x},1-y
 \right]\nonumber\\
 &+ \frac{\g{c_1}\g{b_1-a}\g{1-c_2}\g{1+a+b_2-c_2}}{\g{b_1}\g{1+a-c_2}\g{1+b_2-c_2}\g{c_1-a}}(-x)^{-a} F{}^{2:1;1}_{2:0;0}
  \left[
   \setlength{\arraycolsep}{0pt}
   \begin{array}{c@{{}:{}}c@{;{}}c}
 a,1+a-c_1 & c_2-b_2& b_2\\[1ex]
  c_2,1+a-b_1 & \linefill & \linefill
   \end{array}
   \;\middle|\;
 \frac{1}{x},\frac{1-y}{x}
 \right]
\end{align}
This is sufficient for this series. Now, taking the third series in Eq.\eqref{F_2sol012},
which can also be written as an infinite sum of $_3F_2$ hypergeometric functions:
\begin{multline}
\tilde{F}{}^{2:1;0}_{2:0;0}
  \left[
   \setlength{\arraycolsep}{0pt}
   \begin{array}{c@{{}:{}}c@{;{}}c}
 a+b_2-c_2 ,1+b_2+a-c_1-c_2&b_2,c_2-b_2 & \linefill\\[1ex]
  1+a+b_2-b_1-c_2 ,b_2&\linefill & \linefill
   \end{array}
   \;\middle|\;
 \frac{1-y}{x},x
 \right]=\\\\  \sum_{m=0}^{\infty}\frac{\p{a+b_2-c_2}{m}\p{1+b_2+a-c_2-c_1}{m}}{\p{1+a+b_2-b_1-c_2}{m}m!}\left(\frac{1-y}{x}\right)^{m}\\
 \times
\,_3F_2
\left(
\begin{matrix}
c_2-b_2,\hspace{-0.9cm}& 1-b_2-m,&\hspace{-0.9cm} c_2+b_1-b_2-a-m\\\\
&\hspace{-.3cm}1+c_2-b_2-a-m, & \hspace{-.3cm}c_1+c_2-a-b_2-m   
\end{matrix}
\Bigg {|} {x}
\right) 
\end{multline}
one has, using Eq.\eqref{3F2ac} once more,
\begin{align}
  \tilde{F}{}^{2:1;0}_{2:0;0}
  &\left[
   \setlength{\arraycolsep}{0pt}
   \begin{array}{c@{{}:{}}c@{;{}}c}
 a+b_2-c_2 ,1+b_2+a-c_1-c_2&b_2,c_2-b_2 & \linefill\\[1ex]
  1+a+b_2-b_1-c_2 ,b_2&\linefill & \linefill
   \end{array}
   \;\middle|\;
 \frac{1-y}{x},x
 \right]=\nonumber \\
 & (-x)^{b_2-c_2} \frac{\Gamma \left(1-c_2\right) \Gamma \left(b_1-a\right)  \Gamma \left(-a-b_2+c_2+1\right) \Gamma \left(-a-b_2+c_1+c_2\right)}{\Gamma (1-a) \Gamma \left(1-b_2\right) \Gamma \left(c_1-a\right) \Gamma \left(-a+b_1-b_2+c_2\right)} \nonumber\\
    &\times F{}^{2:1;1}_{2:0;0}
  \left[
   \setlength{\arraycolsep}{0pt}
   \begin{array}{c@{{}:{}}c@{;{}}c}
   a, a-c_1+1&c_2-b_2&b_2\\[1ex]
c_2, a-b_1+1 & \linefill & \linefill
   \end{array}
   \;\middle|\;
 \frac{1}{x},\frac{1-y}{x}
 \right]\nonumber\\
 &+ \frac{\g{b_1+c_2-a-1}\g{1+c_2-b_2-a}\g{c_1+c_2-a-b_2}\g{c_2-1}}{\g{c_2-a}\g{c_2-b_2}\g{c_1+c_2-a-1}\g{b_1+c_2-a-b_2}}(-x)^{-1+b_2}\nonumber\\
 &\times \tilde{F}{}^{1:2;1}_{1:1;0}
  \left[
   \setlength{\arraycolsep}{0pt}
   \begin{array}{c@{{}:{}}c@{;{}}c}
 1-b_2 & 2+a-c_2-c_1,1+a-c_2& b_2\\[1ex]
  2-c_2 & a-c_2-b_1+2 & \linefill
   \end{array}
   \;\middle|\;
 \frac{1}{x},1-y
 \right]\nonumber\\
 &+(-x)^{a-b_1+b_2-c_2}\frac{\Gamma \left(a-b_1\right)  \Gamma \left(a-b_1-c_2+1\right) \Gamma \left(-a-b_2+c_2+1\right) \Gamma \left(-a-b_2+c_1+c_2\right)}{\Gamma \left(1-b_1\right) \Gamma \left(1-b_2\right) \Gamma \left(c_1-b_1\right) \Gamma \left(c_2-b_2\right)} \nonumber\\
 &\times\tilde{F}{}^{1:2;1}_{1:1;0}
  \left[
   \setlength{\arraycolsep}{0pt}
   \begin{array}{c@{{}:{}}c@{;{}}c}
 -a+b_1-b_2+c_2&b_1,b_1-c_1+1&b_2\\[1ex]
  -a+b_1+1& -a+b_1+c_2 & \linefill
   \end{array}
   \;\middle|\;
 \frac{1}{x},1-y
 \right]
\end{align}

\begin{figure}[h!]
  \centering
  \includegraphics[scale=.24]{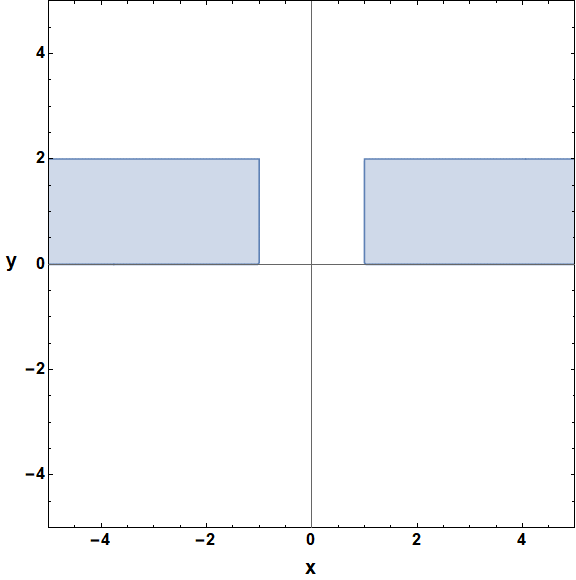}
  \includegraphics[scale=.24]{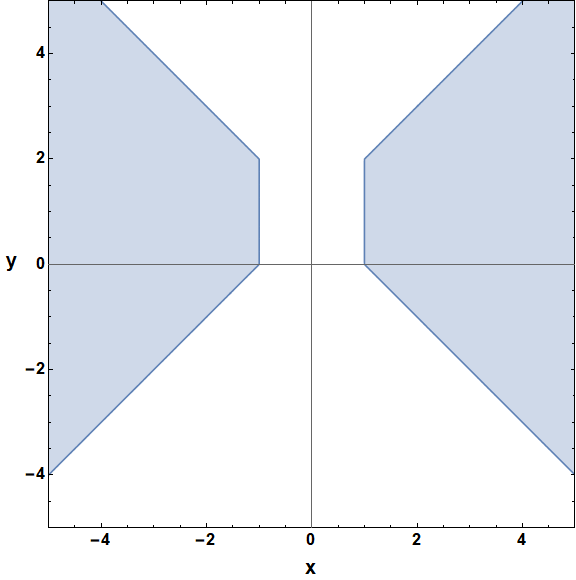}
  \includegraphics[scale=.24]{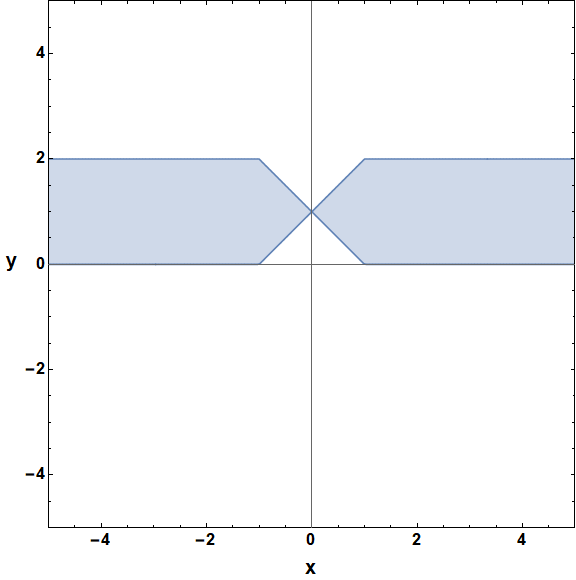}
\caption{\textit{Left:} ROC of $\tilde{F}(...;\frac{1}{x},1-y)$. \textit{Middle:} ROC of $F^{2:1;1}_{2:0;0}(...;\frac{1}{x},\frac{1-y}{x})$. \textit{Right:} ROC of $F^{1:2;1}_{1:1;0}(...;\frac{1-y}{x},1-y)$.\label{F_21inff_t}}
\end{figure}

Substituting the above two results in Eq.\eqref{F_2sol012} and simplifying we get
\begin{align}
&F_2(a,b_1,b_2;c_1,c_2;x,y)=\nonumber\\\nonumber
 &\frac{ \Gamma (c_1) \Gamma (c_2) \Gamma (a-b_1) \Gamma (-a+b_1-b_2+c_2)}{\Gamma (a) \Gamma (c_1-b_1) \Gamma (c_2-b_2) \Gamma (-a+b_1+c_2)}\left(-x\right)^{-b_1}\nonumber\\\nonumber\\\nonumber
&\times \tilde{F}{}^{1:2;1}_{1:1;0}
  \left[
   \setlength{\arraycolsep}{0pt}
   \begin{array}{c@{{}:{}}c@{;{}}c}
 b_1+c_2-a-b_2 & 1+b_1-c_1,b_1& b_2\\[1ex]
  1-a+b_1 & b_1+c_2-a & \linefill
   \end{array}
   \;\middle|\;
 \frac{1}{x},1-y
 \right]\nonumber\\\nonumber\\\nonumber
 &+\frac{\Gamma \left(c_1\right) \Gamma \left(b_1-a\right)}{\Gamma \left(b_1\right) \Gamma \left(c_1-a\right)}\left(-x\right)^{-a}  F^{2:1;1}_{2:0;0}
  \left[
   \setlength{\arraycolsep}{0pt}
   \begin{array}{c@{{}:{}}c@{;{}}c}
 a,1+a-c_1 & c_2-b_2& b_2\\[1ex]
  c_2,1+a-b_1 & \linefill & \linefill
   \end{array}
   \;\middle|\;
 \frac{1}{x},\frac{1-y}{x}
 \right]\nonumber\\\nonumber\\\nonumber
 &+\frac{\Gamma \left(c_1\right) \Gamma \left(c_2\right)  \Gamma \left(a-b_1+b_2-c_2\right)}{\Gamma (a) \Gamma \left(b_2\right) \Gamma \left(c_1-b_1\right)} \left(-x\right)^{-b_1}(1-y)^{-a+b_1-b_2+c_2}\\\nonumber\\
 &\times F{}^{1:2;1}_{1:1;0}
  \left[
   \setlength{\arraycolsep}{0pt}
   \begin{array}{c@{{}:{}}c@{;{}}c}
 c_2-a+b_1 & b_1,1+b_1-c_1& c_2-b_2\\[1ex]
  1+b_1+c_2-b_2-a & c_2-a+b_1 & \linefill
   \end{array}
   \;\middle|\;
 \frac{1-y}{x},1-y
 \right]\label{F_2inf1}
\end{align}
where
$\tilde{F}(...;\frac{1}{x},1-y)$ converges for $|\frac{1}{x}|<1 \wedge |1-y|<1$ (see Fig.\ref{F_21inff_t} \textit{Left}),
$F^{2:1;1}_{2:0;0}(...;\frac{1}{x},\frac{1-y}{x})$ converges for $\Big|\frac{1}{x}\Big|<1 \wedge |\frac{1-y}{x}|<1$ (see Fig.\ref{F_21inff_t} \textit{Middle}) and 
$F^{1:2;1}_{1:1;0}(...;\frac{1-y}{x},1-y)$ converges for $|1-y|<1 \wedge |\frac{1-y}{x}|<1$ (see Fig.\ref{F_21inff_t} \textit{Right}).

It is therefore clear that the analytic continuation in Eq.\eqref{F_2inf1} converges in the region shown in Fig.\ref{F_21inff_t} (\textit{Left}), for real values of $x$ and $y$.

\begin{figure}[h!]
  \centering
 \includegraphics[scale=.24]{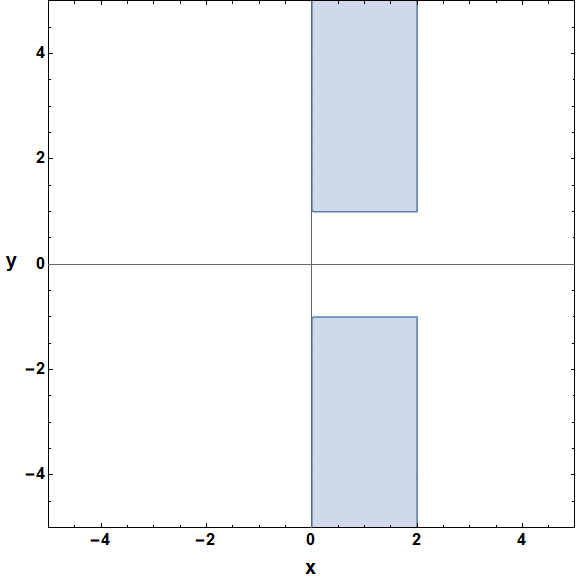}
 \includegraphics[scale=.24]{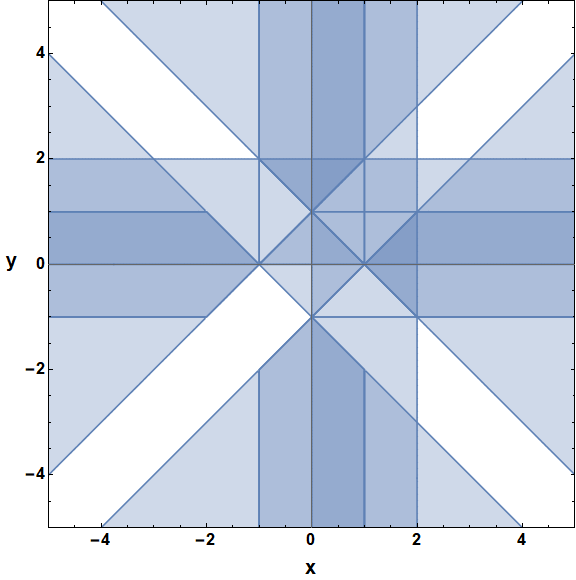}
 \caption{\textit{Left:} See the text. \textit{Right:} Regions covered using all the so far obtained analytic continuations.\label{F_21inff_t2}}
\end{figure}

Using Eq.\eqref{F_2sym} we get another analytic continuation around $(1,\infty)$ which converges in the region $|\frac{1}{y}|<1 \wedge |1-x|<1\wedge |\frac{1-x}{y}|<1$, for real values of $x$ and $y$, see Fig.\ref{F_21inff_t2} (\textit{Left}).\\
\subsection{Using Euler transformations}

In the previous sections, we presented 10 different linear transformations of the Appell $F_2$ double hypergeometric series which analytically continue the latter in various regions of the $(x,y)$ space. 
Now, using these results, it is possible to derive many other analytic continuations. Indeed, if, instead of starting from the series definition of $F_2$, Eq.(\ref{F_2original}), one starts from the alternative series representation given by any of the three Euler transformations shown in Eq.(\ref{F_2_EulerTransforms}), one can use the 10 formulas above to try 30 different ways to derive other analytic continuations. In the total of the various series representations that can be obtained by this exercise (which amounts to 44 if one includes the series definition of $F_2$ and its three Euler transformations), we noticed that 18 are sufficient to cover the whole $(x,y)$ space, except some particular points, namely $(1,0)$, $(0,1)$, $(1,1)$, $(-1,1)$, $(1,-1)$, and $(\frac{1}{2},\frac{1}{2})$. Adding the other analytic continuations does not help to reach the missing points. 
These 18 formulas are listed in the Appendix with plots of their convergence regions, and they form, with the remaining 26 series representations whose expressions are not given explicitly here to lighten the paper, the basis of the \textsc{AppellF2} \textit{Mathematica} package\footnote{The interested reader can find the mathematical expression of each of the 44 series representations of $F_2$ directly from this package.}, presented in Section \ref{package} and dedicated to the numerical evaluation of the Appell $F_2$ function.

\section{The \textsc{\textsc{AppellF2}} \textit{Mathematica} package\label{package}}
Using the series representations of the Appell $F_2$ function derived above and listed in the Appendix, the package \textsc{AppellF2} can find the numerical value of the Appell $F_2$ function, in nonlogarithmic situations (\textit{i.e.} for generic values of the Pochhammer parameters), for arbitrary real values of $x$ and $y$ with the exception of the points that belong to the singular curves of Fig.\ref{RocF2} \textit{Right} for which the result cannot be obtained for all values of the Pochhammer parameters. 

Let us demonstrate the working principle of the package \textsc{AppellF2} below and apply it to some examples later in this section. \textsc{AppellF2} can be used on \textit{Mathematica} v11.3 and beyond.

\subsection{Algorithm of \textsc{AppellF2}}
		The \textsc{AppellF2} package works as follows:
	\begin{enumerate}
		\item Except if the given $(x,y)$ point of evaluation is any of the special points that belong to the singular curves of Fig.\ref{RocF2} \textit{Right}, which are ``Pochhammer dependent", all the series representations of $F_2$ that are valid at the given point are found by the package.  	
		\item Although the same numerical result, for a given precision, will be obtained with any of these series if one sums enough terms, some series will converge faster than others (for instance if the point is not close to the boundary of their convergence region). Therefore, in order to improve the speed of the package, an experimental criterion has been implemented in the code in such a way that the ``best'' series representation, from the convergence point of view, is selected.  The criterion is defined as follows.

A typical series representation of $F_2$ consists of more than one series. 
\begin{align*}
	&F_2(a,b_1,b_2,c_1,c_2,x,y) = \sum_i U_i \\
	&U_i = \sum_{m,n=0}^{\infty}V_i(a,b_1,b_2,c_1,c_2,x,y,m,n)
\end{align*}
For each $U_i$, the package calculates $r_i$ and $s_i$ as written below, for the given values of Pochhammer parameters and $x,y$ (for readibility we have suppressed the dependence of $V_i$ on the Pochhammer parameters and on $x$ and $y$):
\begin{align*}
	s_i = \left|\frac{V_i(m+1,n)}{V_i(m,n)}\right|_{m,n=100}, \hspace{1cm} 	t_i = \left|\frac{V_i(m,n+1)}{V_i(m,n)}\right|_{m,n=100}
\end{align*} 
 The values of $m$ and $n$ have been experimentally chosen to be 100 here but some other values can be used. 
 
The rate of convergence for a series $U_i$ is chosen as $r_i=\sqrt{s_i^2+t_i^2}$. In general, all the series of a given series representation have different rates of convergence. We thus define the rate of convergence of that series representation to be the maximum of the rates of convergence of all of its involved series. Therefore, the rate of convergence of a series representation is
\begin{align*}
R =  \textrm{Max}\{r_i\} 
\end{align*} 

Thus, when comparing the various series representations that are valid at the same $(x,y)$ point, the one that has the smallest $R$ is selected.

\item The numerical evaluation is then performed using partial sums of the best series, for the chosen values of Pochhammer parameters and $x$ and $y$, upto a given number of terms. The output is returned at a given precision.
	\end{enumerate} 
	
\subsection{Demonstration\label{demo}}
We now demonstrate the usage of the package \textsc{AppellF2}. After downloading the package in the same directory as the notebook, it is called as follows

\begin{Verbatim}[fontsize=\small,frame=single]
In[]:SetDirectory[NotebookDirectory[]];
In[]:=<<AppellF2.wl
AppellF2.wl v1.0
Authors : Souvik Bera & Tanay Pathak
\end{Verbatim}

The command \texttt{AppellF2}, computing the numerical value of the $F_2$ function, can be called as,
\begin{Verbatim}[fontsize=\small,frame=single]
In[]:=AppellF2[a, b1, b2, c1, c2, x, y, p, terms, F2show->True]
\end{Verbatim}
Here, \texttt{a,b1,b2,c1,c2} are the Pochhammer values given by the user, \texttt{x,y} is the point of evaluation, \texttt{p} is the required precision of the output (it gives the number of desired significant digits) and \texttt{terms} is the number of terms in the numerical summation for each summation index.  \texttt{F2show} is an option with default value \texttt{False}. When it is made \texttt{True}, one can see the evaluation of the summation in real time. 

As an example, the below command finds the value of $F_2$ upto $4$ significant digits at the point $(-2.311, 5.322)$ with Pochhammer values $a=2.2345, b_1=3.363, b_2=0.242, c_1=8.3452$ and $c_2=0.657$ by evaluating the summation upto $100$ terms for both the summation indices.
\begin{Verbatim}[fontsize=\small,frame=single]
In[]:=AppellF2[2.2345, 3.363, 0.242, 8.3452, 0.657,-2.311, 5.322, 4, 100, F2show->False]
\end{Verbatim}
For this call, the package gives
\begin{Verbatim}[fontsize=\small,frame=single]
valid series:{{10},{15},{18},{26},{29},{43}}
convergence rates:{{0.59,10},{0.66,26},{0.68,18},{0.94,29},{0.95,43},{1.09,15}}
selected series: 10
	
Out[]=0.09334 - 0.06847 I
\end{Verbatim}
One can see that there are 6 series representations valid at the point $(-2.311, 5.322)$ and that series $\#10$ is the best converging series that is chosen for the evaluation.
 
One will note that the specific command \texttt{F2findall} can be directly called to find which series are valid at a given point. For the example above,
\begin{Verbatim}[fontsize=\small,frame=single]
In[]:=F2findall[{-2.311, 5.322}]
	
Out[]={10, 15, 18, 26, 29, 43}
\end{Verbatim}

In order to see the expression of any of the 44 series representations used in the package, and its region of convergence, the command \texttt{F2expose} can be used.

\begin{Verbatim}[fontsize=\small,frame=single]
In[]:=F2expose[15]
	
Out[]= {Abs[(-1 + x + y)/x]< 1&&Abs[x/(-1 + x)]<1,(1/(m! n!))((1-x)^{-a}) Gamma[c2]...}
\end{Verbatim}
The output above, which, due to its length, has been partly suppressed here, is a list containing the ROC of the series representation followed by the expression of the corresponding series (here $\#15$ which, from the correspondence Table \ref{denomination} of the appendix, is $S_7$, see Eq.(\ref{S7})).

\begin{figure}[h]
  \centering
  \includegraphics[scale=.5]{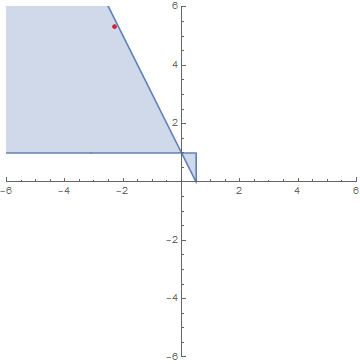}
  \caption{Output of \texttt{F2ROC[\string{-2.311, 5.322\string}, 15, \{-6, 6\}]}. The red dot is $(x,y)$=$(-2.311,5.322)$.}\label{F2figdemo}
\end{figure}

The command \texttt{F2ROC} gives the plot of the ROC of the series representation \texttt{\#}, along with the point \texttt{(x,y)}, for a given \texttt{range}. 
\begin{Verbatim}[fontsize=\small,frame=single]
F2ROC[{x,y}, #, range]
\end{Verbatim}
For instance, the call \texttt{F2ROC[\string{-2.311, 5.322\string}, 15, \string{-6, 6\string}]}
gives the output shown in Fig. \ref{F2figdemo}.

Finally, the user can chose which series representation he wants for the evaluation at a given point, using the command \texttt{F2evaluate}. For example,
\begin{Verbatim}[fontsize=\small,frame=single]
In[]:=F2evaluate[10, {2.2345, 3.363, 0.242, 8.3452, 0.657, -2.311, 5.322}, 10, 100]
	
Out[]=0.09333639793 - 0.06847416686 I
\end{Verbatim}
gives the values of $F_2(2.2345, 3.363, 0.242; 8.3452, 0.657; -2.311, 5.322)$ derived from the series representation $\#10$ at 10 precision with 100 terms for each summation index.

\subsection{Numerical tests\label{numerics}}

We have tested the \textsc{AppellF2} package by computing 200 randomly generated points: 100 points with the ranges of random values of the Pochhammer parameters $a, b_1, b_2, c_1$ and $c_2$ and of the $x$ and $y$ variables being $[-7,7]$, and 100 points with random complex values of the Pochhammer parameters (with real and imaginary parts in the range $[-7,7]$).

For each of these 200 points, the series representations that are valid there all match numerically. As some of the series converge faster than others at a given point, it is sometimes necessary to increase the number of terms in the partial sums for those series that converge slowly. 

We have also compared these results to the \textit{Maple} inbuilt \textit{AppellF2} function. There is a very good agreement for 168 points. In the 32 remaining points there are 16 points for which \textit{Maple} does not give a numerical result, and 16 points for which there is a discrepancy between \textit{Maple}'s results and those obtained from the series of \textsc{AppellF2} that are valid at those points. The problematic points with real values of the Pochhammer parameters are shown in Table \ref{MapleVSAppellF2} (we do not show the points with complex Pochhammer parameters for lack of space). It is obvious that for points $\#2, \#8$ and $\#17$ in this table, \textit{Maple} gives incorrect results, as the numerical evaluation of $F_2$ at these points should be real, whereas \textit{Maple} gives a nonzero imaginary part. Indeed, points $\#2, \#8$ and $\#17$ lie in the convergence region of the third Euler transformation of $F_2$ (Eq.(\ref{S3})) which asks for an evaluation of $F_2$ in the convergence region of its usual series definition, which has to be real. Furthermore, the prefactor of this Euler transformation is not evaluated on its branch cut for the corresponding $x$ and $y$ values and therefore does not yield an imaginary part. 
This analysis is corroborated for point $\#2$ by the numerical evaluation of the Euler integral representation of $F_2$ given in Eq.(\ref{F_2_Euler2}), which confirms the result obtained from \textsc{AppellF2} (for points $\#8$ and $\#17$ the integral does not converge).

Further investigation is needed to better understand the discrepancies for the other 6 points of Table \ref{MapleVSAppellF2}. 

\begin{longtable}{|c|c|c|c|c|}
\hline
Serial  & Pochhammer parameters & \textit{Maple} output & \textsc{AppellF2} output & Series  \\
no & and $x,y$& &  & number\\\hline
\endfirsthead
\endhead
1 & \begin{tabular}[c]{@{}c@{}}$a$= -4.49158729455734\\ $b_1$= 4.69491717746224\\ $b_2$= -2.67898515537678\\ $c_1$= 2.54939072003598\\ $c_2$= 1.89372308769086\\ $x$= -0.657865707164980\\ $y$ = 1.11972469394233\end{tabular} & $183.83 - 0.00072 i$ & $183.83$ & \begin{tabular}[c]{@{}c@{}}40\\ 32\\ 4\\ 15\\ 16\end{tabular} \\ \hline
2 &\begin{tabular}[c]{@{}c@{}}$a$= -5.87056003391116\\ $b_1$= 4.33993527730256\\ $b_2$= 1.44218908732163\\ $c_1$= 3.12652020729955\\ $c_2$=1.52984418542146\\ $x$= -6.55177221618387\\ $y$= -6.79935054310963\end{tabular} & $1.171 \times 10^7 + 0.019 i$ & $1.171 \times 10^7$ & \begin{tabular}[c]{@{}c@{}}34\\ 26\\ 16\\ 35\\36\end{tabular} \\ \hline
3&\begin{tabular}[c]{@{}c@{}}$a$= 3.72583256450429\\  $b_1$= 2.11602255447865\\ $b_2$= -3.02238392715598\\ $c_1$= -4.73946336645648\\ $c_2$= 6.30095725032474\\  $x$= -2.59888480330968\\ $y$= 3.17343904351674\end{tabular} & - & $-0.46 - 0.018 i$ & \begin{tabular}[c]{@{}c@{}}29\\ 15\\ 43\\ 16\\ 28\end{tabular} \\ \hline
4&\begin{tabular}[c]{@{}c@{}}$a$= -2.88936562201761\\  $b_1$= 3.45488861254925\\  $b_2$= -5.90441801674065\\  $c_1$= 2.41900748028973\\ $c_2$= -4.22539504741494\\  $x$= 5.41515683798479\\ $y$= -3.96256437728474\end{tabular} & - & $-0.25 - 0.055 i$ & \begin{tabular}[c]{@{}c@{}}44\\ 11\\ 27\\ 17\\ 16\\ 28\end{tabular} \\ \hline
5&\begin{tabular}[c]{@{}c@{}}$a$= 3.91112960454197\\  $b_1$= 0.943419234377764\\ $b_2$= -0.939723628477704\\  $c_1$= -6.05906265997090\\  $c_2$= -6.10615861109657\\  $x$= -2.12080578449250\\ $y$= 2.86626245522156\end{tabular} & - & $5.99 \times 10^{20} + 4.67\times 10^{18} i$ & \begin{tabular}[c]{@{}c@{}}29\\ 15\\ 43\\ 16\\ 28\end{tabular} \\ \hline
6&\begin{tabular}[c]{@{}c@{}}$a$= 6.31031029746575\\ $b_1$= -4.81608319937911\\ $b_2$= 0.134215670834948\\  $c_1$= -2.75396269319432\\  $c_2$= -3.75042638259086\\  $x$= -4.08388732316032\\ $y$= 1.81702135884373\end{tabular} & $2.03 \times 10^{13} + 6.45\times 10^{13} i$ & $3.17 \times10^{12} - 1.10 \times 10^{11} i$ & \begin{tabular}[c]{@{}c@{}}33\\ 16\\ 38\\ 11\\ 6\\ 15\end{tabular} \\ \hline
7&\begin{tabular}[c]{@{}c@{}}$a$= 4.14277514262421\\ $b_1$= -6.43436403118499\\ $b_2$= 1.60386793716277\\ $c_1$= -6.87730424656771\\ $c_2$= -5.67535554477487\\ $x$= -0.646110168140300\\ $y$= 0.817740014591525\end{tabular} & $4.05\times 10^{19} +  5.93 \times10^{12} i$ & $4.05\times 10^{19}$ & \begin{tabular}[c]{@{}c@{}}12\\ 4\\ 13\\ 38 \\ 14\end{tabular} \\ \hline
8&\begin{tabular}[c]{@{}c@{}}$a$= 3.35171139159466\\ $b_1$= -0.509725596574174\\ $b_2$= -0.913836915342344\\ $c_1$= -3.32588271257136\\ $c_2$= 0.168816510623319\\ $x$= -2.29531801533183\\ $y$ = -6.06415712186627\end{tabular} & $-61.36 + 0.01 i$ & $-61.38$ & \begin{tabular}[c]{@{}c@{}}26\\ 10\\ 18\\ 34\\ 35\\ 19\\ 36\end{tabular} \\ \hline
9&\begin{tabular}[c]{@{}c@{}}$a$= -5.01240784115629\\ $b_1$= -4.94200818581766\\ $b_2$= 6.99477562102917\\ $c_1$= 6.65313744284692\\ $c_2$= -1.96099117581162\\ $x$= 2.92126097205082\\ $y$= -1.31245113310376\end{tabular} & $6.00\times 10^6 - 0.00032 i$ & $6.00\times 10^6$ & \begin{tabular}[c]{@{}c@{}}44\\ 27\\ 17\\ 28\\ 11\\ 16\end{tabular} \\ \hline
10&\begin{tabular}[c]{@{}c@{}}$a$= -0.981118466281753\\  $b_1$= 4.55280800772390\\ $b_2$= 1.43404196228123\\ $c_1$= 2.84087159624645\\  $c_2$= 6.00107528411102\\ $x$= -5.36758744763326\\  $y$= 6.61806381273987\end{tabular} & - & $7.72 + 0.000027 i$ & \begin{tabular}[c]{@{}c@{}}43\\ 15\\ 29\\ 10\\ 26\\ 18\end{tabular} \\ \hline
11&\begin{tabular}[c]{@{}c@{}}$a$= 1.04079628966533\\  $b_1$= 0.999310189508378\\  $b_2$= 3.59329558885096\\ $c_1$= -6.25832679000047\\ $c_2$= -4.02905455852754\\ $x$= 1.81023829524087\\ $y$= 2.00777521482951\end{tabular} & $-7.72\times 10^{16} + 1.23\times 10^{18} i$ & $-122497.46 - 140113.63 i$ & \begin{tabular}[c]{@{}c@{}}37\\ 38\\ 7\\ 28\end{tabular} \\ \hline
12&\begin{tabular}[c]{@{}c@{}}$a$= -3.26985266196408\\  $b_1$= 0.380743118208180\\ $b_2$= 2.02474976684470\\ $c_1$= -1.31514385444273\\  $c_2$= .83951473440144\\ $x$= 5.35725173812456\\  $y$= -4.07499617412362\end{tabular} & - & $304.17 + 38.11 i$ & \begin{tabular}[c]{@{}c@{}}44\\ 27\\ 17\\ 11\\ 16\\28\end{tabular} \\ \hline
13&\begin{tabular}[c]{@{}c@{}}$a$= -6.17654955276504\\ $b_1$= 3.21556912170448\\  $b_2$= -2.26411156484076\\ $c_1$= 5.34290089035759\\  $c_2$= -2.32233932454304\\  $x$= 3.61217179206409\\ $y$= -2.53790197142651\end{tabular} & - & $158.84 + 62.05 i$ & \begin{tabular}[c]{@{}c@{}}44\\ 27\\ 17\\ 11\\ 28\\ 16\end{tabular} \\ \hline
14&\begin{tabular}[c]{@{}c@{}}$a$= 2.31197860013321\\  $b_1$= -0.666975465151342\\  $b_2$= -5.24476192259412\\  $c_1$= -3.16508771091695\\  $c_2$= 6.55592102157901\\  $x$= 6.31953155096413\\ $y$= -6.36985521062664\end{tabular} & - & $7322.40 - 12654.38 i$ & \begin{tabular}[c]{@{}c@{}}18\\ 26\\ 17\\ 27\end{tabular} \\ \hline
15&\begin{tabular}[c]{@{}c@{}}$a$= 1.16647934045583\\  $b_1$= -2.56252103706461\\  $b_2$= 5.52623207935986\\  $c_1$= 6.58905357119552\\  $c_2$= 5.53232577263389\\ $x$= -2.05715925643916\\ $y$= 3.15407000645672\end{tabular} & - & $-0.23 + 0.13 i$ & \begin{tabular}[c]{@{}c@{}}43\\ 15\\ 29\\ 10\\ 18\\26\end{tabular} \\ \hline
16&\begin{tabular}[c]{@{}c@{}}$a$= 4.26170736723804\\ $b_1$= 2.41512776824820\\ $b_2$= -3.60520211982802\\  $c_1$= -2.44037707125234\\ $c_2$= -3.72147640617149\\ $x$= -2.08902304321602\\ $y$= 4.40568570030866\end{tabular} & $-7.46\times 10^{13} + 7.54\times 10^{13} i$ & $ -13729.68 + 34149.43 i$ & \begin{tabular}[c]{@{}c@{}}10\\ 43\\ 18\\ 26\\ 29\\ 15\end{tabular} \\ \hline
17&\begin{tabular}[c]{@{}c@{}}$a$= -3.36021432698409\\  $b_1$= 6.63749440272489\\  $b_2$= -6.58339249087694\\ $c_1$= -2.02579013838810\\  $c_2$= 6.18081281041145\\  $x$= -4.71272838790961\\  $y$= -6.11479355971970\end{tabular} & $6.03 \times 10^9 + 2.74 \times 10^9 i$ & $-3.20\times 10^6$ & \begin{tabular}[c]{@{}c@{}}34\\ 26\\ 10\\ 18\\ 35\\ 36\end{tabular} \\ \hline
\caption{In the set of 100 randomly generated real points that we have carefully tested, these are those for which there is a discrepancy between \textit{Maple} and \textsc{AppellF2} (or no result from \textit{Maple}). Each example, in the table above, is described by a set of values for $(a,b_1,b_2,c_1,c_2,x,y)$. The results from the package at these points are obtained from partial sums with 300 terms for each summation index $m$ and $n$, with precision of 15 significant digits (the results are however truncated with less significant digits in the table to improve its readability). It should also be noted that the order in which the series in the last column are given corresponds to the decreasing order of convergence rate.
\label{MapleVSAppellF2}}
\end{longtable}

\section{Conclusions}

The Appell $F_2$ is  an important hypergeometric function of two variables, which can be linked to 10 of the 14 complete double hypergeometric functions of order two.

In this work, we have carried out a comprehensive analysis of this function and built its implementation for \textit{Mathematica} in the form of the \textsc{AppellF2} package presented in Section \ref{package} and provided as an ancillary file to this paper.
Our method starts from the original series definition of $F_2(a,b_1,b_2;c_1,c_2;x,y)$, which has a limited range of validity, on which we apply the transformation theory following Olsson's approach \cite{Olsson-64}. In this way, we derived a set of 43 linear transformations for $F_2$.
These formulas, which are valid for generic (\textit{i.e.} for nonlogarithmic cases) complex values of the $a, b_1, b_2, c_1, c_2$ parameters, can collectively cover the entire real $(x, y)$ space, as concluded from the study of their regions of convergence, except on a few particular points.

In fact, 18 formulas in this set of 44 series representations of $F_2$ are sufficient to cover the real $(x,y)$ space, but we have incorporated all the 44 in the \textsc{AppellF2} package in order to improve its convergence efficiency.
We have also carefully studied the behavior of these formulas on their branch cuts and we have given their expressions there in a consistent way.

The usage of the package has been explained in Section \ref{package} where the numerical checks have also been carried out, as described in Section \ref{numerics}, to confirm the consistency of our results, internally and also by a comparison with the existing \textit{AppellF2} inbuilt function of \textit{Maple} with which we find disagreement at several instances, as shown in Table \ref{MapleVSAppellF2}.

Let us note here that the \textsc{AppellF2} package can be used to develop a \textit{Mathematica} realization of the  ten second order complete hypergeometric functions in two variables which $F_2$ is linked to, which is a work in progress \cite{workinprogress}. As another extension of the present study, we will also consider, in a subsequent work, the logarithmic situation where some of the Pochhammer parameters can be identical. 

\section*{Acknowledgments}
 S.F. thanks Centre for High Energy Physics, Indian Institutes of Science of Bangalore, where this work was initiated, for hospitality.

\section*{Appendix: Series representations of the Appell $F_2$ function}
	
We list in this Appendix the 18 series representations $S_i, (i=1,...,18)$ of $F_2$ that can collectively cover the $(x,y)$ real space and we plot their regions of convergence. The remaining 26 series used in the package can be obtained from the latter by calling the \texttt{F2expose[]} command, as explained in Section \ref{demo}. We refer the reader to Eqs.(\ref{KdFnotation}) and (\ref{MirrorKdFnotation}) for the notation of Kamp\'e de F\'eriet series and their mirror partners used in the expressions presented in this appendix.	

In the following 18 formulas (whose corresponding denomination in the \textsc{AppellF2} package are listed in Table \ref{denomination}), which can also be seen as functional relations, the involved series appear in general multiplied by prefactors which have the form of powers of linear rational functions of $x$ and $y$ whose exponent can be fractional. It may happen that these prefactors, if not carefully considered when evaluated on their branch cuts, lead to inconsistencies between the different formulas that are valid at the same points. Therefore, to obtain a matching it is necessary to proceed to the rewriting of some of the prefactors (as this has been performed for the Gauss $_2F_1$ case in \cite{B&S}) in a way which is equivalent everywhere except on the branch cuts where the rewriting gives the expected behavior, in agreement with the conventions of \textit{Mathematica}.

To perform this rewriting of the prefactors, we have defined a conditional function, denoted by brackets, as follows:
\begin{align}
        \left\langle \left(f(x,y)\right)^a \right\rangle  =
        \left\{
\begin{array}{rcr}
\left(f(x,y)\right)^a & \text{if} & \text{Condition} \\
\left(\frac{1}{f(x,y)}\right)^{-a} & \text{else} &  
\end{array}
\right.
\end{align}
where $a$ is any linear combination of Pochhammer parameters. 
The conditions for all the prefactors that appear in the expressions of 44 series representations used in the \textsc{AppellF2} package are summarized in Table \ref{conditions}:

\begin{table}[h]
\begin{center}
        \begin{tabular}{ |c| c| }
              \hline  Argument of the prefactor & Condition \\ \hline 
                $\frac{x}{y-1}$, $\frac{x+y-1}{y-1}$ & $x-y+1>0 $\\ 
                $-\frac{y}{x+y-1}$, $\frac{x-1}{x+y-1}$ & $x-y-1>0 $\\ 
                $\frac{y}{x-1}$, $\frac{x+y-1}{x-1}$ & $-x+y+1>0  $\\ 
                $-\frac{x}{x+y-1}$, $\frac{y-1}{x+y-1}$ & $-x+y-1>0  $  \\
                $\frac{1}{1-x}$, $\frac{1}{1-y}$, $-\frac{x}{x-1}$,  $-\frac{y}{y-1}$ &$ \text{False}$\\
              \hline
        \end{tabular}
                              \caption{Conditions of rewriting of the possible prefactors that appear in the 44 series representations of $F_2$ used in the \textsc{AppellF2} package.\label{conditions}}
                              \end{center}
\end{table}
As an example,
\begin{align}
        \left\langle \left(\frac{x}{y-1}\right)^{-a-b_2+c_2} \right\rangle  =
        \left\{
\begin{array}{rcr}
\left(\frac{x}{y-1}\right)^{-a-b_2+c_2} & \text{if} & x-y+1>0 \\
\left(\frac{y-1}{x}\right)^{a+b_2-c_2} & \text{else} &  
\end{array}
\right.
\end{align}

\begin{table}
\begin{tabular}{|c|c|}
 \hline
    Series representations in the appendix &  Series representations $\#$ in the \textsc{AppellF2} package  \\
    \hline
     $S_1$ & 1 \\
    $S_2$ & 23\\
    $S_3$ & 34\\
    \hline
    $S_4$ & 14\\
    $S_5$ & 25\\
    $S_6$ & 4\\
     \hline
    $S_7$ & 15\\
    $S_8$ & 37\\
    $S_9$ & 5\\
     \hline
    $S_{10}$ & 27\\
    $S_{11}$ & 38\\
    $S_{12}$ & 6\\ 
    \hline
    $S_{13}$ & 17\\
    $S_{14}$ & 7\\
    $S_{15}$ & 29\\
     \hline
    $S_{16}$ & 40\\ 
    $S_{17}$ & 8\\
    $S_{18}$ & 9\\
     \hline
\end{tabular}
\caption{Denomination of the 18 series of the appendix in the \textsc{AppellF2} package.\label{denomination}}
\end{table}
	\subsection*{Series representation $S_1$}
	$S_1$ is the original series
	\begin{align}
		S_1=F_2(a,b_1,b_2;c_1,c_2;x,y)=\sum_{m=0}^{\infty}\sum_{n=0}^{\infty}\frac{(a)_{m+n} (b_1)_m (b_2)_n }{(c_1)_m (c_2)_n m! n! }x^m y^n
	\end{align}
	whose region of convergence is $| x|+| y|<1$ (see Fig. \ref{ROCS1} \textit{left}).

	\subsection*{Series representation $S_2$}
	$S_2$ is one of the Euler transformations of $F_2$.
	\begin{align}
		S_2&=(1 - y)^{-a}F_2\left(a,b_1,c_2-b_2;c_1,c_2;\frac{x}{1-y},\frac{y}{y-1}\right)
	\end{align}
Region of convergence: $| \frac{x}{1-y}| +| \frac{y}{y-1}| <1$ (see Fig. \ref{ROCS1} \textit{middle}).

	\subsection*{Series representation $S_3$}
\begin{figure}[h]
		\centering
		\includegraphics[scale=.24]{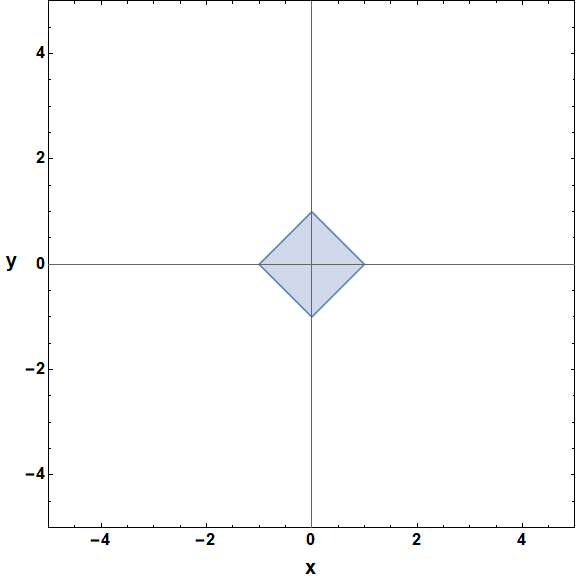}
		\includegraphics[scale=.24]{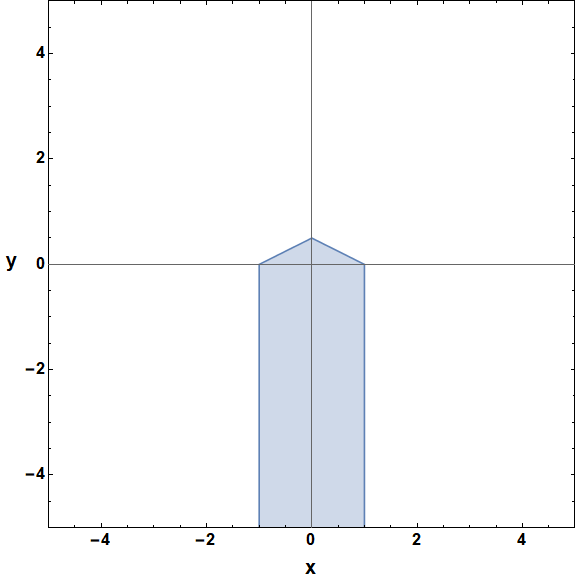}
		\includegraphics[scale=.24]{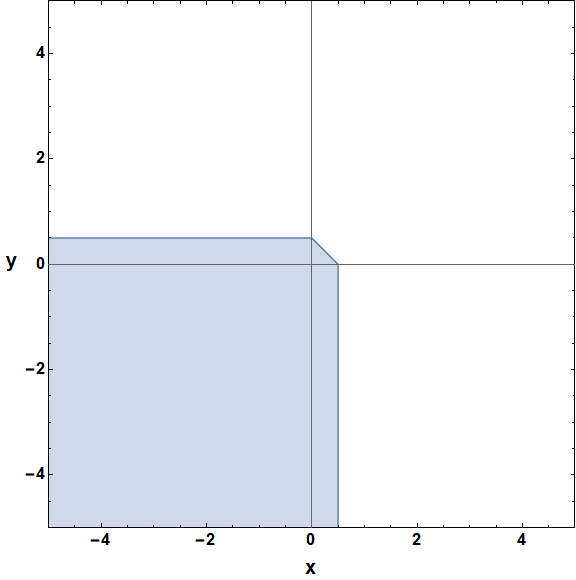}
		\caption{Convergence regions of $S_1$ (\textit{left}), $S_2$ (\textit{middle}) and $S_3$ (\textit{right}) for real values of $x$ and $y$.\label{ROCS1}}
	\end{figure}

		$S_3$ is another Euler transformation of $F_2$.
	\begin{align}
		S_3&=(1-x-y)^{-a}F_2\left(a,c_1-b_1,c_2-b_2;c_1,c_2;\frac{x}{x+y-1},\frac{y}{x+y-1}\right)\label{S3}
	\end{align}
Region of convergence: $| \frac{x}{x+y-1}| +| \frac{y}{x+y-1}| <1$ (see Fig. \ref{ROCS1} \textit{right}).
	
	\subsection*{Series representation $S_4$}

	$S_4$ is obtained by applying Eq.(\ref{F_2sol101}) on the first Euler transformation of $F_2$ (first line of Eq.(\ref{F_2_EulerTransforms}), \textit{i.e.} $S_2$).
	\begin{align}
		S_4&= (1-x)^{-a}\left\langle\left(\frac{1}{1-x}\right)^{-a+b_1}\right\rangle \frac{\Gamma (c_1) \Gamma (a-b_1)}{\Gamma (a) \Gamma (c_1-b_1)}\tilde{F}{}^{1:1;2}_{1:0;1}
  \left[
   \setlength{\arraycolsep}{0pt}
   \begin{array}{c@{{}:{}}c@{;{}}c}
  c_1-a & b_1 & b_2,a-c_1+1\\[1ex]
   b_1-a+1 & - & c_2
   \end{array}
   \;\middle|\;
 \frac{1}{1-x},y
 \right]\nonumber\\
		&+(1-x)^{-a} \frac{\Gamma (c_1) \Gamma (b_1-a)}{\Gamma (b_1) \Gamma (c_1-a)}{F}{}^{1:1;2}_{1:0;1}
  \left[
   \setlength{\arraycolsep}{0pt}
   \begin{array}{c@{{}:{}}c@{;{}}c}
  a & c_1-b_1 & b_2,a-c_1+1\\[1ex]
   a-b_1+1 & - & c_2
   \end{array}
   \;\middle|\;
 \frac{1}{1-x},\frac{y}{1-x}
 \right]
	\end{align}
	Region of convergence: $| 1-x| >1\land | y| <1$ (see Fig. \ref{ROCS2} \textit{left}).

	\subsection*{Series representation $S_5$}

	$S_5$ is obtained by applying Eq.(\ref{F_2sol101}) on the second Euler transformation of $F_2$ (second line of Eq.(\ref{F_2_EulerTransforms})).
	\begin{align}
		S_5&=(1-y)^{-a} \frac{\Gamma (c_1) \Gamma (-a-b_1+c_1)}{\Gamma (c_1-a) \Gamma (c_1-b_1)} {F}{}^{1:1;2}_{1:0;1}
  \left[
   \setlength{\arraycolsep}{0pt}
   \begin{array}{c@{{}:{}}c@{;{}}c}
  a & b_1 & c_2-b_2,a-c_1+1\\[1ex]
   a+b_1-c_1+1 & - & c_2
   \end{array}
   \;\middle|\;
 \frac{x+y-1}{y-1},\frac{y}{y-1}
 \right]\nonumber\\
		&+ (1-y)^{-a}\left\langle\left(\frac{x+y-1}{y-1}\right)^{-a-b_1+c_1}\right\rangle\frac{ \Gamma (c_1) \Gamma (a+b_1-c_1) }{\Gamma (a) \Gamma (b_1)}\nonumber\\
		&\times\tilde{F}{}^{1:1;2}_{1:0;1}
  \left[
   \setlength{\arraycolsep}{0pt}
   \begin{array}{c@{{}:{}}c@{;{}}c}
  c_1-a & c_1-b_1 & c_2-b_2,a-c_1+1\\[1ex]
   -a-b_1+c_1+1 & - & c_2
   \end{array}
   \;\middle|\;
 \frac{x+y-1}{y-1},\frac{y}{x+y-1}
 \right]
	\end{align}
	Region of convergence: $| \frac{y}{x+y-1}| <1\land | \frac{x+y-1}{y-1}| <1$ (see Fig. \ref{ROCS2} \textit{middle}).

	\subsection*{Series representation $S_6$}
\begin{figure}[h]
		\centering
		\includegraphics[scale=.24]{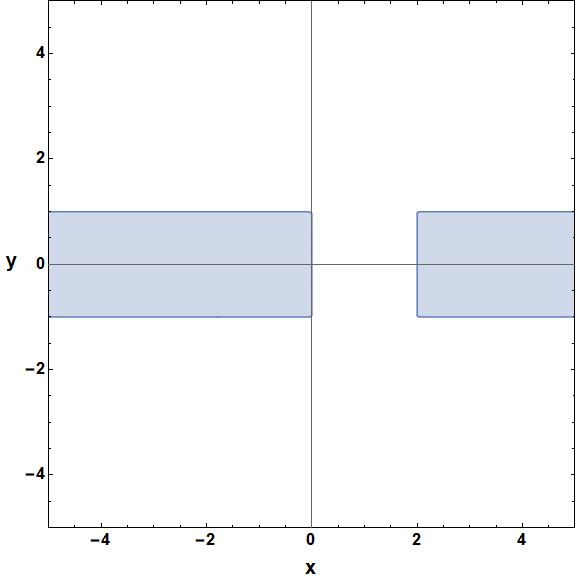}
		\includegraphics[scale=.24]{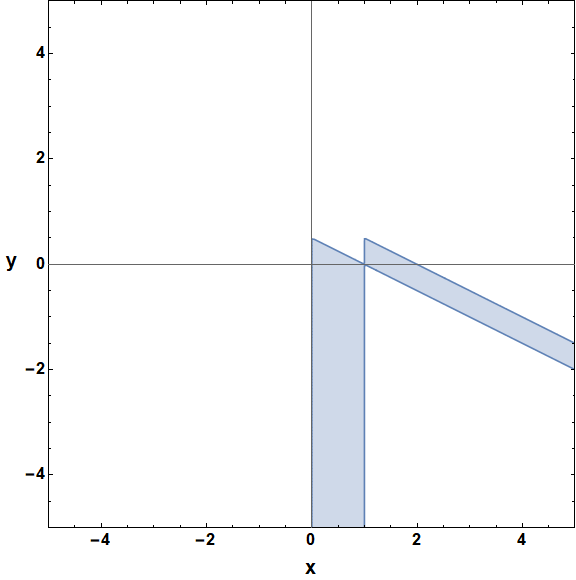}
		\includegraphics[scale=.24]{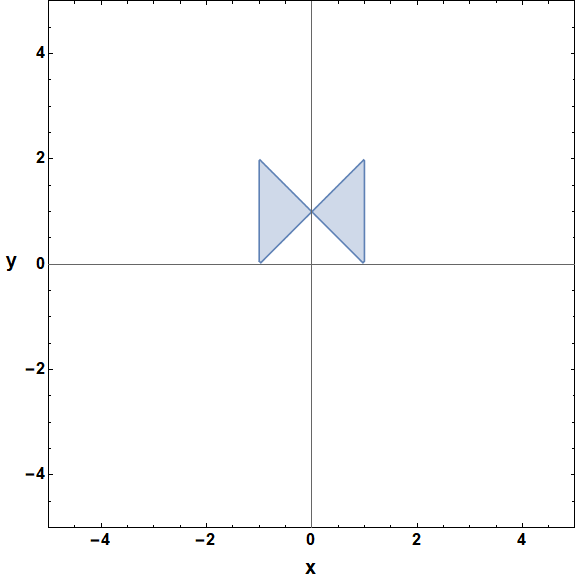}
		\caption{Convergence regions of $S_4$ (\textit{left}), $S_5$ (\textit{middle}) and $S_6$ (\textit{right}) for real values of $x$ and $y$.\label{ROCS2}}
	\end{figure}

	$S_6$ is Eq.(\ref{F_2sol012}).
	\begin{align}
		S_6&= \frac{\Gamma (c_2) \Gamma (-a-b_2+c_2)}{\Gamma (c_2-a) \Gamma (c_2-b_2)} {F}{}^{1:2;1}_{1:1;0}
  \left[
   \setlength{\arraycolsep}{0pt}
   \begin{array}{c@{{}:{}}c@{;{}}c}
  a & a-c_2+1,b_1 & b_2\\[1ex]
   a+b_2-c_2+1 & c_1 & -
   \end{array}
   \;\middle|\;
 x,1-y
 \right]\nonumber\\
		&+(1-y)^{-a-b_2+c_2}\left\langle\left(\frac{x}{y-1}\right)^{-a-b_2+c_2} \right\rangle\frac{\Gamma (c_1) \Gamma (c_2)  \Gamma (a+b_2-c_2) \Gamma (-a+b_1-b_2+c_2) }{\Gamma (a) \Gamma (b_1) \Gamma (b_2) \Gamma (-a-b_2+c_1+c_2)}\nonumber\\
		&\times\tilde{F}{}^{2:1;1}_{2:0;0}
  \left[
   \setlength{\arraycolsep}{0pt}
   \begin{array}{c@{{}:{}}c@{;{}}c}
  a+b_2-c_2, a+b_2-c_2-c_1+1 & b_2 & c_2-b_2\\[1ex]
   b_2,a-b_1+b_2-c_2+1 & - & -
   \end{array}
   \;\middle|\;
 \frac{1-y}{x},x
 \right]\nonumber\\
		&+\left\langle\left(\frac{x}{y-1}\right)^{-b_1}\right\rangle (1-y)^{-a-b_2+c_2}\frac{\Gamma (c_1) \Gamma (c_2)  \Gamma (a-b_1+b_2-c_2)}{\Gamma (a) \Gamma (b_2) \Gamma (c_1-b_1)}\nonumber\\
		&\times{F}{}^{1:2;1}_{1:1;0}
  \left[
   \setlength{\arraycolsep}{0pt}
   \begin{array}{c@{{}:{}}c@{;{}}c}
  -a+b_1+c_2 & b_1,b_1-c_1+1 & c_2-b_2\\[1ex]
   -a+b_1-b_2+c_2+1 & -a+b_1+c_2 & -
   \end{array}
   \;\middle|\;
 \frac{1-y}{x},1-y
 \right]
	\end{align}
	Region of convergence: $| \frac{1-y}{x}| <1\land | x| <1$ (see Fig. \ref{ROCS2} \textit{right}).
	
	\subsection*{Series representation $S_7$}
	$S_7$ is obtained by applying Eq.(\ref{F_2sol012}), \textit{i.e.} $S_6$, on the first Euler transformation of $F_2$ (first line of Eq.(\ref{F_2_EulerTransforms}), \textit{i.e.} $S_2$).
	\begin{align}
		S_7&=(1-x)^{-a}\frac{\Gamma (c_2) \Gamma (-a-b_2+c_2)}{\Gamma (c_2-a) \Gamma (c_2-b_2)} {F}{}^{1:2;1}_{1:1;0}
  \left[
   \setlength{\arraycolsep}{0pt}
   \begin{array}{c@{{}:{}}c@{;{}}c}
  a & a-c_2+1, c_1-b_1& b_2\\[1ex]
   a+b_2-c_2+1 & c_1 & - 
   \end{array}
   \;\middle|\;
 \frac{x}{x-1},\frac{x+y-1}{x-1}
 \right]\nonumber\\
		&+ (1-x)^{-a}\left\langle\left(-\frac{x}{x+y-1}\right)^{-a-b_2+c_2}\right\rangle\left\langle\left(\frac{x+y-1}{x-1}\right)^{-a-b_2+c_2}\right\rangle \frac{\Gamma (c_1) \Gamma (c_2) \Gamma (a+b_2-c_2)}{\Gamma (a) \Gamma (b_2) \Gamma (c_1-b_1) } \nonumber\\ &\times\frac{ \Gamma (-a-b_1-b_2+c_1+c_2)  }{\Gamma (-a-b_2+c_1+c_2)}\tilde{F}{}^{2:1;1}_{2:0;0}
  \left[
   \setlength{\arraycolsep}{0pt}
   \begin{array}{c@{{}:{}}c@{;{}}c}
  a+b_2-c_2, a+b_2-c_1-c_2+1  & b_2 & c_2-b_2\\[1ex]
   a+b_1+b_2-c_1-c_2+1, b_2 & - & -
   \end{array}
   \;\middle|\;
 \frac{x+y-1}{x},\frac{x}{x-1}
 \right]\nonumber\\
		&+ (1-x)^{-a} \left\langle\left(\frac{x+y-1}{x-1}\right)^{-a-b_2+c_2}\right\rangle\left\langle\left(-\frac{x}{x+y-1}\right)^{b_1-c_1}\right\rangle\frac{ \Gamma (c_1) \Gamma (c_2) \Gamma (a+b_1+b_2-c_1-c_2) }{\Gamma (a) \Gamma (b_1) \Gamma (b_2)}\nonumber\\
		&\times{F}{}^{1:2;1}_{1:1;0}
  \left[
   \setlength{\arraycolsep}{0pt}
   \begin{array}{c@{{}:{}}c@{;{}}c}
  -a-b_1+c_1+c_2 & 1-b_1,c_1-b_1 & c_2-b_2\\[1ex]
   -a-b_1-b_2+c_1+c_2+1 & -a-b_1+c_1+c_2 & -
   \end{array}
   \;\middle|\;
 \frac{x+y-1}{x},\frac{x+y-1}{x-1}
 \right]\label{S7}
	\end{align}
	Region of convergence: $| \frac{x+y-1}{x}| <1\land | \frac{x}{x-1}| <1$ (see Fig. \ref{ROCS3} \textit{left}).

	\subsection*{Series representation $S_8$}
	$S_8$ is obtained by applying Eq.(\ref{F_2sol012}), \textit{i.e.} $S_6$, on the third Euler transformation of $F_2$ (third line of Eq.(\ref{F_2_EulerTransforms}), \textit{i.e.} $S_3$).
	\begin{align}
		S_8&=(-x-y+1)^{-a}\left\langle\left(-\frac{x}{x-1}\right)^{b_2-a}\right\rangle\left\langle\left(\frac{x-1}{x+y-1}\right)^{b_2-a}\right\rangle\frac{\Gamma (c_1) \Gamma (c_2)  \Gamma (a-b_2)  \Gamma (-a-b_1+b_2+c_1) }{\Gamma (a) \Gamma (c_1-b_1) \Gamma (c_2-b_2) \Gamma (-a+b_2+c_1)} \nonumber\\
		&\times\tilde{F}{}^{2:1;1}_{2:0;0}
  \left[
   \setlength{\arraycolsep}{0pt}
   \begin{array}{c@{{}:{}}c@{;{}}c}
  a-b_2, a-b_2-c_1+1 & c_2-b_2 & b_2 \\[1ex]
   c_2-b_2,a+b_1-b_2-c_1+1 & - & -
   \end{array}
   \;\middle|\;
 \frac{x-1}{x},\frac{x}{x+y-1}
 \right]\nonumber\\
		&+ (-x-y+1)^{-a}\frac{\Gamma (c_2) \Gamma (b_2-a) }{\Gamma (b_2) \Gamma (c_2-a)}{F}{}^{1:2;1}_{1:1;0}
  \left[
   \setlength{\arraycolsep}{0pt}
   \begin{array}{c@{{}:{}}c@{;{}}c}
  a & a-c_2+1, c_1-b_1 & c_2-b_2\\[1ex]
   a-b_2+1 & c_1 & -
   \end{array}
   \;\middle|\;
 \frac{x}{x+y-1},\frac{x-1}{x+y-1}
 \right]\nonumber\\
		&+  (-x-y+1)^{-a}\left\langle\left(-\frac{x}{x-1}\right)^{b_1-c_1}\right\rangle\left\langle\left(\frac{x-1}{x+y-1}\right)^{b_2-a}\right\rangle\frac{\Gamma (c_1) \Gamma (c_2)  \Gamma (a+b_1-b_2-c_1) }{\Gamma (a) \Gamma (b_1) \Gamma (c_2-b_2)} \nonumber\\
		&\times{F}{}^{1:2;1}_{1:1;0}
  \left[
   \setlength{\arraycolsep}{0pt}
   \begin{array}{c@{{}:{}}c@{;{}}c}
  -a-b_1+c_1+c_2 & 1-b_1,c_1-b_1 & b_2\\[1ex]
   -a-b_1+b_2+c_1+1 & -a-b_1+c_1+c_2 & -
   \end{array}
   \;\middle|\;
\frac{x-1}{x},\frac{x-1}{x+y-1}
 \right]\label{S8}
	\end{align}
	Region of convergence: $| \frac{x-1}{x}| <1\land | \frac{x}{x+y-1}| <1$ (see Fig. \ref{ROCS1} \textit{middle}).

\subsection*{Series representation $S_9$}
\begin{figure}[h]
		\centering
		\includegraphics[scale=.24]{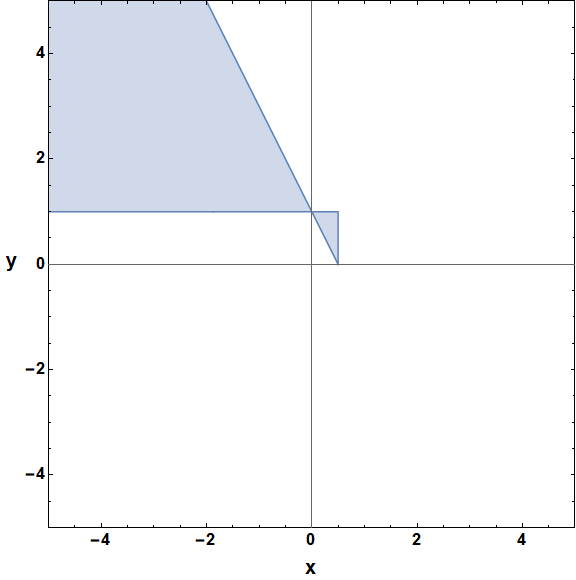}
		\includegraphics[scale=.24]{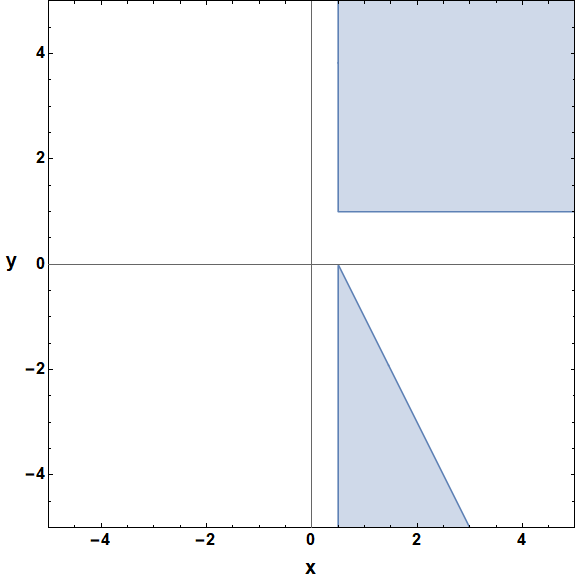}
		\includegraphics[scale=.24]{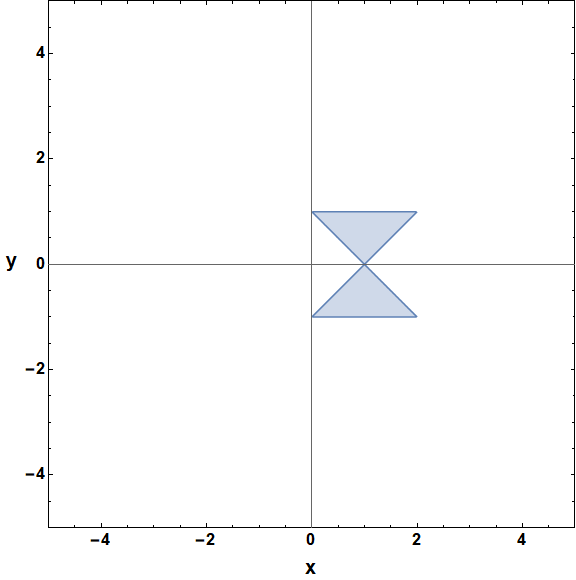}
		\caption{Convergence regions of $S_7$ (\textit{left}), $S_8$ (\textit{middle}) and $S_9$ (\textit{right}) for real values of $x$ and $y$.\label{ROCS3}}
	\end{figure}

$S_9$ is the symmetrical partner of Eq.(\ref{F_2sol012}) obtained from Eq.(\ref{F_2sym}).
	\begin{align}
		S_9&=(1-x)^{-a-b_1+c_1}\left\langle\left(\frac{y}{x-1}\right)^{-b_2}\right\rangle\frac{\Gamma (c_1) \Gamma (c_2)   \Gamma (a+b_1-b_2-c_1)}{\Gamma (a) \Gamma (b_1) \Gamma (c_2-b_2)}\nonumber\\
		&\times{F}{}^{1:1;2}_{1:0;1}
  \left[
   \setlength{\arraycolsep}{0pt}
   \begin{array}{c@{{}:{}}c@{;{}}c}
  -a+b_2+c_1 & c_1-b_1 & b_2 , b_2-c_2+1\\[1ex]
   -a-b_1+b_2+c_1+1 & - & -a+b_2+c_1
   \end{array}
   \;\middle|\;
 1-x,\frac{1-x}{y}
 \right]\nonumber\\
		&+ (1-x)^{-a-b_1+c_1}\left\langle\left(\frac{y}{x-1}\right)^{-a-b_1+c_1}\right\rangle\frac{\Gamma (c_1) \Gamma (c_2)  \Gamma (a+b_1-c_1) \Gamma (-a-b_1+b_2+c_1) }{\Gamma (a) \Gamma (b_1) \Gamma (b_2) \Gamma (-a-b_1+c_1+c_2)}\nonumber\\
		&\times\tilde{F}{}^{2:1;1}_{2:0;0}
  \left[
   \setlength{\arraycolsep}{0pt}
   \begin{array}{c@{{}:{}}c@{;{}}c}
  1-b_1, -a-b_1+b_2+c_1 & c_1-b_1 & b_1\\[1ex]
   -a-b_1+c_1+1, -a-b_1+c_1+c_2 & - & -
   \end{array}
   \;\middle|\;
 y,\frac{1-x}{y}
 \right]\nonumber\\
		&+ \frac{\Gamma (c_1) \Gamma (-a-b_1+c_1)}{\Gamma (c_1-a) \Gamma (c_1-b_1)}{F}{}^{1:1;2}_{1:0;1}
  \left[
   \setlength{\arraycolsep}{0pt}
   \begin{array}{c@{{}:{}}c@{;{}}c}
  a & b_1 & b_2,a-c_1+1\\[1ex]
   a+b_1-c_1+1 & - & c_2
   \end{array}
   \;\middle|\;
 1-x,y
 \right]
	\end{align}
	Region of convergence: $| \frac{1-x}{y}| <1\land | y| <1$ (see Fig. \ref{ROCS3} \textit{right}).
	
	\subsection*{Series representation $S_{10}$}
	$S_{10}$ is obtained by applying the symmetrical partner Eq.(\ref{F_2sol012}), \textit{i.e.} $S_9$, on the second Euler transformation of $F_2$ (second line of Eq.(\ref{F_2_EulerTransforms})).
	\begin{align}
		S_{10}&=(1-y)^{-a}\left\langle\left(-\frac{y}{x+y-1}\right)^{-a-b_1+c_1}\right\rangle\left\langle \left(\frac{x+y-1}{y-1}\right)^{-a-b_1+c_1}\right\rangle \frac{\Gamma (c_1) \Gamma (c_2) \Gamma (a+b_1-c_1)  }{\Gamma (a) \Gamma (b_1) \Gamma (c_2-b_2) }\nonumber\\ 
		&\times\frac{ \Gamma (-a-b_1-b_2+c_1+c_2) }{ \Gamma (-a-b_1+c_1+c_2)}\tilde{F}{}^{2:1;1}_{2:0;0}
  \left[
   \setlength{\arraycolsep}{0pt}
   \begin{array}{c@{{}:{}}c@{;{}}c}
  1-b_1, -a-b_1-b_2+c_1+c_2 & c_1-b_1 & b_1\\[1ex]
   -a-b_1+c_1+1, -a-b_1+c_1+c_2 & - & -
   \end{array}
   \;\middle|\;
\frac{y}{y-1},\frac{x+y-1}{y}
 \right]\nonumber\\
		&+ (1-y)^{-a}\frac{ \Gamma (c_1) \Gamma (-a-b_1+c_1)}{\Gamma (c_1-a) \Gamma (c_1-b_1)}{F}{}^{1:1;2}_{1:0;1}
  \left[
   \setlength{\arraycolsep}{0pt}
   \begin{array}{c@{{}:{}}c@{;{}}c}
  a & b_1 & a-c_1+1, c_2-b_2\\[1ex]
   a+b_1-c_1+1 & - & c_2
   \end{array}
   \;\middle|\;
\frac{x+y-1}{y-1},\frac{y}{y-1}
 \right]\nonumber\\
		&+(1-y)^{-a}\left\langle\left(\frac{x+y-1}{y-1}\right)^{-a-b_1+c_1}\right\rangle\left\langle\left(-\frac{y}{x+y-1}\right)^{b_2-c_2}\right\rangle \frac{ \Gamma (c_1) \Gamma (c_2) \Gamma (a+b_1+b_2-c_1-c_2)  }{\Gamma (a) \Gamma (b_1) \Gamma (b_2)}\nonumber\\
		&\times{F}{}^{1:1;2}_{1:0;1}
  \left[
   \setlength{\arraycolsep}{0pt}
   \begin{array}{c@{{}:{}}c@{;{}}c}
  -a-b_2+c_1+c_2 & c_1-b_1 & 1-b_2, c_2-b_2\\[1ex]
   -a-b_1-b_2+c_1+c_2+1 & - & -a-b_2+c_1+c_2
   \end{array}
   \;\middle|\;
\frac{x+y-1}{y-1},\frac{x+y-1}{y}
 \right]
	\end{align}
	Region of convergence: $| \frac{x+y-1}{y}| <1\land | \frac{y}{y-1}| <1$ (see Fig. \ref{ROCS4} \textit{left}).

	\subsection*{Series representation $S_{11}$}
	$S_{11}$ is obtained by applying the symmetrical partner Eq.(\ref{F_2sol012}), \textit{i.e.} $S_9$, on the third Euler transformation of $F_2$ (third line of Eq.(\ref{F_2_EulerTransforms}), \textit{i.e.} $S_3$).
	\begin{align}
		S_{11}&= (-x-y+1)^{-a}\left\langle\left(\frac{y}{1-y}\right)^{b_1-a}\right\rangle\left\langle \left(\frac{y-1}{x+y-1}\right)^{b_1-a}\right\rangle\frac{\Gamma (c_1) \Gamma (c_2) \Gamma (a-b_1)  \Gamma (-a+b_1-b_2+c_2)}{\Gamma (a) \Gamma (c_1-b_1) \Gamma (c_2-b_2) \Gamma (-a+b_1+c_2)}\nonumber\\
		&\times\tilde{F}{}^{2:1;1}_{2:0;0}
  \left[
   \setlength{\arraycolsep}{0pt}
   \begin{array}{c@{{}:{}}c@{;{}}c}
  b_1-c_1+1, -a+b_1-b_2+c_2 & b_1 & c_1-b_1\\[1ex]
   -a+b_1+1, -a+b_1+c_2 & - & -
   \end{array}
   \;\middle|\;
\frac{y}{x+y-1},\frac{y-1}{y}
 \right]\nonumber\\
		&+ (-x-y+1)^{-a}\frac{\Gamma (c_1) \Gamma (b_1-a) }{\Gamma (b_1) \Gamma (c_1-a)}{F}{}^{1:1;2}_{1:0;1}
  \left[
   \setlength{\arraycolsep}{0pt}
   \begin{array}{c@{{}:{}}c@{;{}}c}
  a & c_1-b_1 & a-c_1+1, c_2-b_2\\[1ex]
   a-b_1+1 & - & c_2
   \end{array}
   \;\middle|\;
\frac{y-1}{x+y-1},\frac{y}{x+y-1}
 \right]\nonumber\\
		&+(-x-y+1)^{-a} \left\langle\left(\frac{y}{1-y}\right)^{b_2-c_2}\right\rangle\left\langle\left(\frac{y-1}{x+y-1}\right)^{b_1-a}\right\rangle \frac{\Gamma (c_1) \Gamma (c_2)  \Gamma (a-b_1+b_2-c_2) }{\Gamma (a) \Gamma (b_2) \Gamma (c_1-b_1)}\nonumber\\
		& \times{F}{}^{1:1;2}_{1:0;1}
  \left[
   \setlength{\arraycolsep}{0pt}
   \begin{array}{c@{{}:{}}c@{;{}}c}
  -a-b_2+c_1+c_2 & b_1 & 1-b_2, c_2-b_2\\[1ex]
   -a+b_1-b_2+c_2+1 & - & -a-b_2+c_1+c_2
   \end{array}
   \;\middle|\;
\frac{y-1}{x+y-1},\frac{y-1}{y}
 \right]
	\end{align}
	Region of convergence: $| \frac{y-1}{y}| <1\land | \frac{y}{x+y-1}| <1$ (see Fig. \ref{ROCS4} \textit{middle}).

	\subsection*{Series representation $S_{12}$}
\begin{figure}[h]
		\centering
		\includegraphics[scale=.24]{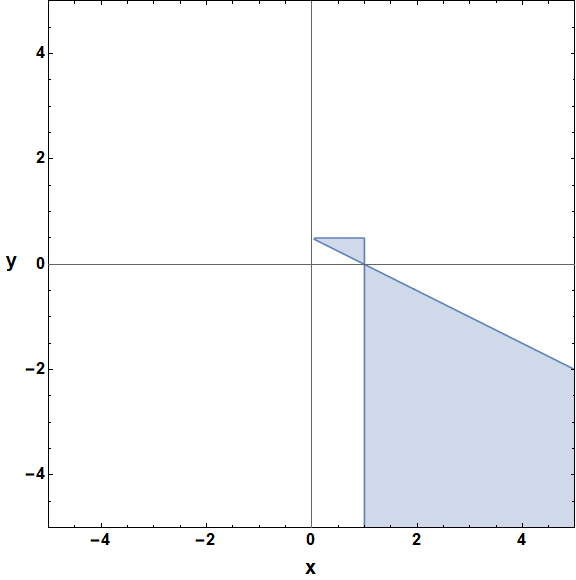}
		\includegraphics[scale=.24]{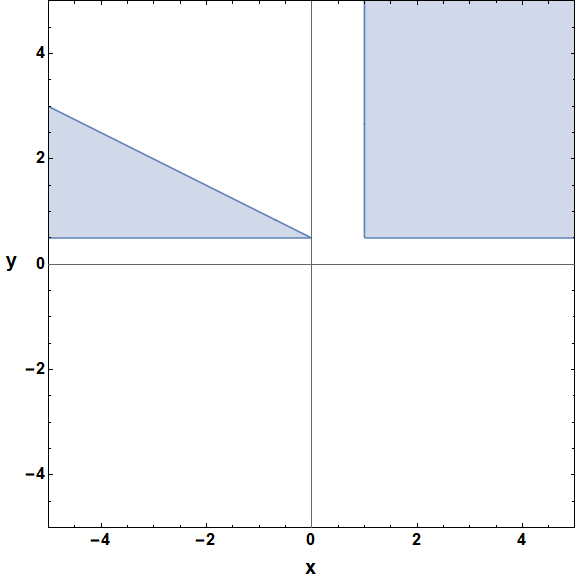}
		\includegraphics[scale=.24]{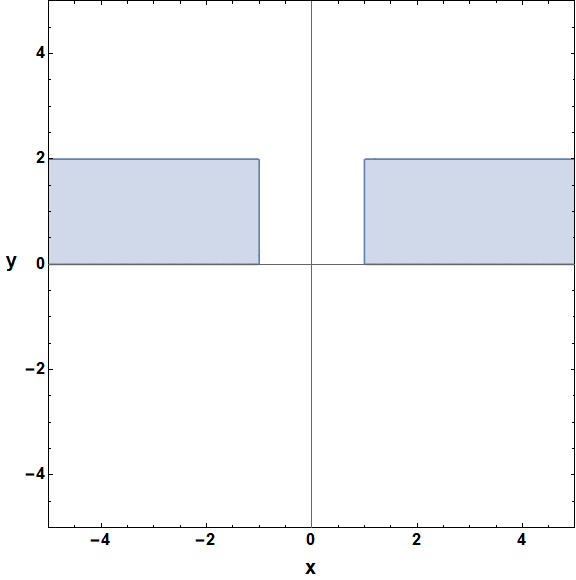}
		\caption{Convergence regions of $S_{10}$ (\textit{left}), $S_{11}$ (\textit{middle}) and $S_{12}$ (\textit{right}) for real values of $x$ and $y$.\label{ROCS4}}
	\end{figure}

	$S_{12}$ is Eq.(\ref{F_2inf1}). 
	\begin{align}
		S_{12}&=(-x)^{-a}\frac{ \Gamma (c_1) \Gamma (b_1-a)}{\Gamma (b_1) \Gamma (c_1-a)} {F}{}^{2:1;1}_{2:0;0}
  \left[
   \setlength{\arraycolsep}{0pt}
   \begin{array}{c@{{}:{}}c@{;{}}c}
  a, a-c_1+1 & c_2-b_2 & b_2\\[1ex]
   c_2, a-b_1+1 & - & -
   \end{array}
   \;\middle|\;
\frac{1}{x},\frac{1-y}{x}
 \right] \nonumber\\
		&+(-x)^{-b_1} \frac{ \Gamma (c_1) \Gamma (c_2) \Gamma (a-b_1) \Gamma (-a+b_1-b_2+c_2)}{\Gamma (a) \Gamma (c_1-b_1) \Gamma (c_2-b_2) \Gamma (-a+b_1+c_2)}\tilde{F}{}^{1:2;1}_{1:1;0}
  \left[
   \setlength{\arraycolsep}{0pt}
   \begin{array}{c@{{}:{}}c@{;{}}c}
  -a+b_1-b_2+c_2 & b_1, b_1-c_1+1 & b_2\\[1ex]
   -a+b_1+1 & -a+b_1+c_2 & -
   \end{array}
   \;\middle|\;
\frac{1}{x},1-y
 \right]\nonumber\\
		&+ (-x)^{-b_1}(1-y)^{-a+b_1-b_2+c_2} \frac{ \Gamma (c_1) \Gamma (c_2) \Gamma (a-b_1+b_2-c_2)}{\Gamma (a) \Gamma (b_2) \Gamma (c_1-b_1)}\nonumber\\
		&\times{F}{}^{1:2;1}_{1:1;0}
  \left[
   \setlength{\arraycolsep}{0pt}
   \begin{array}{c@{{}:{}}c@{;{}}c}
  -a+b_1+c_2 & b_1, b_1-c_1+1 & c_2-b_2\\[1ex]
   -a+b_1-b_2+c_2+1 & -a+b_1+c_2 & -
   \end{array}
   \;\middle|\;
\frac{1-y}{x},1-y
 \right]
	\end{align}
	Region of convergence: $| \frac{1}{x}| <1\land | 1-y| <1$ (see Fig. \ref{ROCS4} \textit{right}).
	
	\subsection*{Series representation $S_{13}$}
	$S_{13}$ is obtained by applying Eq.(\ref{F_2inf1}), \textit{i.e. $S_{13}$}, on the first Euler transformation of $F_2$ (first line of Eq.(\ref{F_2_EulerTransforms}), \textit{i.e.} $S_2$). 
	\begin{align}
		S_{13}&=(1-x)^{-a}\left\langle\left(-\frac{x}{x-1}\right)^{-a}\right\rangle \frac{ \Gamma (c_1) \Gamma (-a-b_1+c_1)}{\Gamma (c_1-a) \Gamma (c_1-b_1)}{F}{}^{2:1;1}_{2:0;0}
  \left[
   \setlength{\arraycolsep}{0pt}
   \begin{array}{c@{{}:{}}c@{;{}}c}
  a, a-c_1+1 & c_2-b_2 & b_2\\[1ex]
   c_2, a+b_1-c_1+1 & - & -
   \end{array}
   \;\middle|\;
\frac{x-1}{x},\frac{x+y-1}{x}
 \right]\nonumber\\
		&+ (1-x)^{-a}\left\langle\left(-\frac{x}{x-1}\right)^{b_1-c_1}\right\rangle\frac{ \Gamma (c_1) \Gamma (c_2)  \Gamma (a+b_1-c_1) \Gamma (-a-b_1-b_2+c_1+c_2)}{\Gamma (a) \Gamma (b_1) \Gamma (c_2-b_2) \Gamma (-a-b_1+c_1+c_2)}\nonumber\\
		& \times\tilde{F}{}^{1:2;1}_{1:1;0}
  \left[
   \setlength{\arraycolsep}{0pt}
   \begin{array}{c@{{}:{}}c@{;{}}c}
  -a-b_1-b_2+c_1+c_2 & 1-b_1, c_1-b_1 & b_2\\[1ex]
   -a-b_1+c_1+1 & -a-b_1+c_1+c_2 & -
   \end{array}
   \;\middle|\;
\frac{x-1}{x},\frac{x+y-1}{x-1}
 \right]\nonumber\\
		&+ (1-x)^{-a} \left\langle\left(-\frac{x}{x-1}\right)^{b_1-c_1}\right\rangle\left\langle\left(\frac{x+y-1}{x-1}\right)^{-a-b_1-b_2+c_1+c_2}\right\rangle\frac{ \Gamma (c_1) \Gamma (c_2) \Gamma (a+b_1+b_2-c_1-c_2) }{\Gamma (a) \Gamma (b_1) \Gamma (b_2)}\nonumber\\
		& \times{F}{}^{1:2;1}_{1:1;0}
  \left[
   \setlength{\arraycolsep}{0pt}
   \begin{array}{c@{{}:{}}c@{;{}}c}
  -a-b_1+c_1+c_2 & 1-b_1, c_1-b_1 & c_2-b_2\\[1ex]
   -a-b_1-b_2+c_1+c_2+1 & -a-b_1+c_1+c_2 & -
   \end{array}
   \;\middle|\;
\frac{x+y-1}{x},\frac{x+y-1}{x-1}
 \right]
	\end{align}
	Region of convergence: $| \frac{x-1}{x}| <1\land | \frac{x+y-1}{x-1}| <1$ (see Fig. \ref{ROCS5} \textit{left}).

	\subsection*{Series representation $S_{14}$}
	$S_{14}$  is the symmetrical partner of Eq.(\ref{F_2inf1}) obtained from Eq.(\ref{F_2sym}).
	\begin{align}
		S_{14}&=(-y)^{-a}\frac{ \Gamma (c_2) \Gamma (b_2-a)}{\Gamma (b_2) \Gamma (c_2-a)} {F}{}^{2:1;1}_{2:0;0}
  \left[
   \setlength{\arraycolsep}{0pt}
   \begin{array}{c@{{}:{}}c@{;{}}c}
  a, a-c_2+1 & b_1 & c_1-b_1\\[1ex]
  c_1, a-b_2+1 & - & -
   \end{array}
   \;\middle|\;
\frac{1-x}{y},\frac{1}{y}
 \right]\nonumber\\
		&+ (-y)^{-b_2}\frac{ \Gamma (c_1) \Gamma (c_2) \Gamma (a-b_2) \Gamma (-a-b_1+b_2+c_1)}{\Gamma (a) \Gamma (c_1-b_1) \Gamma (c_2-b_2) \Gamma (-a+b_2+c_1)}\tilde{F}{}^{1:1;2}_{1:0;1}
  \left[
   \setlength{\arraycolsep}{0pt}
   \begin{array}{c@{{}:{}}c@{;{}}c}
  a-b_2 & b_1 & b_2, b_2-c_2+1\\[1ex]
   a+b_1-b_2-c_1+1 & - & -a+b_2+c_1
   \end{array}
   \;\middle|\;
1-x,\frac{1}{y}
 \right]\nonumber\\
		&+(-y)^{-b_2}(1-x)^{-a-b_1+b_2+c_1} \frac{ \Gamma (c_1) \Gamma (c_2)  \Gamma (a+b_1-b_2-c_1)}{\Gamma (a) \Gamma (b_1) \Gamma (c_2-b_2)}\nonumber\\
		& \times{F}{}^{1:1;2}_{1:0;1}
  \left[
   \setlength{\arraycolsep}{0pt}
   \begin{array}{c@{{}:{}}c@{;{}}c}
  -a+b_2+c_1 & c_1-b_1 & b_2, b_2-c_2+1\\[1ex]
   -a-b_1+b_2+c_1+1 & - & -a+b_2+c_1
   \end{array}
   \;\middle|\;
1-x,\frac{1-x}{y}
 \right]
	\end{align}
	Region of convergence: $| 1-x| <1\land | y| >1$ (see Fig. \ref{ROCS5} \textit{middle}).

	\subsection*{Series representation $S_{15}$}
\begin{figure}[h]
		\centering
		\includegraphics[scale=.24]{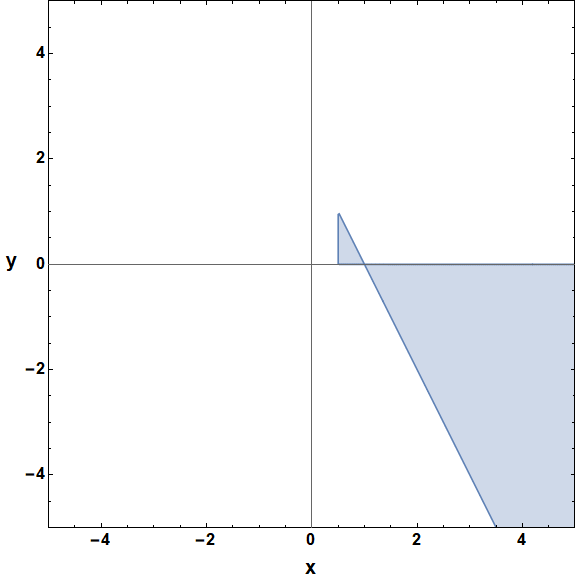}
		\includegraphics[scale=.24]{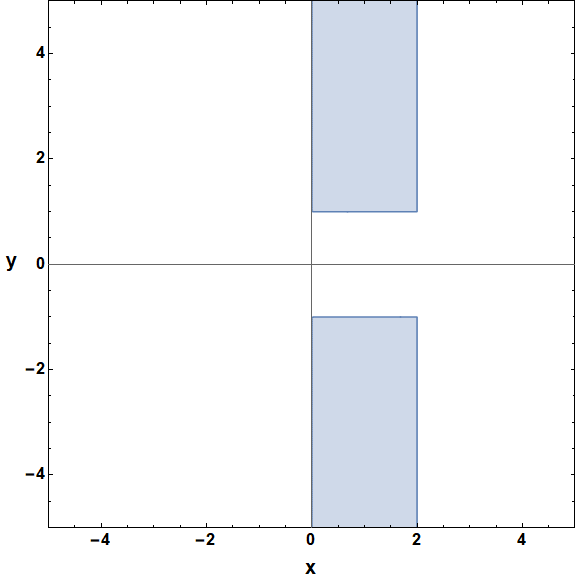}
		\includegraphics[scale=.24]{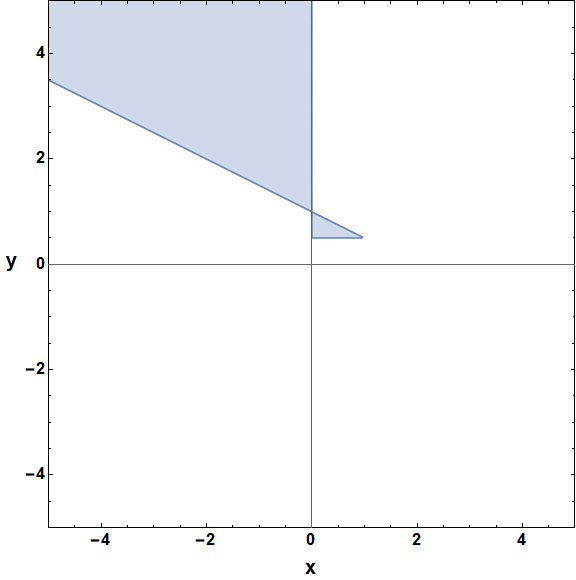}
		\caption{Convergence regions of $S_{13}$ (\textit{left}), $S_{14}$ (\textit{middle}) and $S_{15}$ (\textit{right}) for real values of $x$ and $y$.\label{ROCS5}}
	\end{figure}

	$S_{15}$ is obtained by applying the symmetrical partner of Eq.(\ref{F_2inf1}), \textit{i.e.} $S_{14}$, on the second Euler transformation of $F_2$ (second line of Eq.(\ref{F_2_EulerTransforms})).
	\begin{align}
		S_{15}&=(1-y)^{-a}\left\langle\left(-\frac{y}{y-1}\right)^{-a} \right\rangle\frac{ \Gamma (c_2) \Gamma (-a-b_2+c_2)}{\Gamma (c_2-a) \Gamma (c_2-b_2)}{F}{}^{2:1;1}_{2:0;0}
  \left[
   \setlength{\arraycolsep}{0pt}
   \begin{array}{c@{{}:{}}c@{;{}}c}
  a, a-c_2+1 & b_1 & c_1-b_1\\[1ex]
  c_1, a+b_2-c_2+1 & - & -
   \end{array}
   \;\middle|\;
\frac{x+y-1}{y},\frac{y-1}{y}
 \right] \nonumber\\
		&+ (1-y)^{-a} \left\langle\left(-\frac{y}{y-1}\right)^{b_2-c_2}\right\rangle \frac{\Gamma (c_1) \Gamma (c_2) \Gamma (a+b_2-c_2) \Gamma (-a-b_1-b_2+c_1+c_2) }{\Gamma (a) \Gamma (b_2) \Gamma (c_1-b_1) \Gamma (-a-b_2+c_1+c_2)}\nonumber\\ 
		& \times\tilde{F}{}^{1:1;2}_{1:0;1}
  \left[
   \setlength{\arraycolsep}{0pt}
   \begin{array}{c@{{}:{}}c@{;{}}c}
  a+b_2-c_2 & b_1 & 1-b_2,c_2-b_2\\[1ex]
  a+b_1+b_2-c_1-c_2+1 & - & -a-b_2+c_1+c_2
   \end{array}
   \;\middle|\;
\frac{x+y-1}{y-1},\frac{y-1}{y}
 \right] \nonumber\\
		&+ (1-y)^{-a}\left\langle\left(-\frac{y}{y-1}\right)^{b_2-c_2}\right\rangle \left\langle\left(\frac{y-1}{x+y-1}\right)^{a+b_1+b_2-c_1-c_2}\right\rangle\frac{ \Gamma (c_1) \Gamma (c_2)  \Gamma (a+b_1+b_2-c_1-c_2)}{\Gamma (a) \Gamma (b_1) \Gamma (b_2)} \nonumber\\
		&\times{F}{}^{1:1;2}_{1:0;1}
  \left[
   \setlength{\arraycolsep}{0pt}
   \begin{array}{c@{{}:{}}c@{;{}}c}
  -a-b_2+c_1+c_2 & c_1-b_1 & 1-b_2,c_2-b_2\\[1ex]
  -a-b_1-b_2+c_1+c_2+1 & - & -a-b_2+c_1+c_2
   \end{array}
   \;\middle|\;
\frac{x+y-1}{y-1},\frac{x+y-1}{y}
 \right]
	\end{align}
Region of convergence: $| \frac{y-1}{y}| <1\land | \frac{x+y-1}{y-1}| <1$ (see Fig. \ref{ROCS5} \textit{right}).

	\subsection*{Series representation $S_{16}$}
	$S_{16}$ is obtained by applying the symmetrical partner of Eq.(\ref{F_2inf1}), \textit{i.e.} $S_{14}$, on the third Euler transformation of $F_2$ (third line of Eq.(\ref{F_2_EulerTransforms}), \textit{i.e.} $S_3$).
\begin{align}
		S_{16}&=(-x-y+1)^{-a}\left\langle\left(-\frac{y}{x+y-1}\right)^{-a}\right\rangle\frac{\Gamma (c_2)   \Gamma (-a-b_2+c_2)}{\Gamma (c_2-a) \Gamma (c_2-b_2)}\nonumber\\
		& \times{F}{}^{2:1;1}_{2:0;0}
  \left[
   \setlength{\arraycolsep}{0pt}
   \begin{array}{c@{{}:{}}c@{;{}}c}
  a, a-c_2+1 & c_1-b_1 & b_1\\[1ex]
  c_1, a+b_2-c_2+1 & - & -
   \end{array}
   \;\middle|\;
\frac{y-1}{y},\frac{x+y-1}{y}
 \right]\nonumber\\
		&+(-x-y+1)^{-a}\left\langle\left(-\frac{y}{x+y-1}\right)^{b_2-c_2}\right\rangle \frac{\Gamma (c_1) \Gamma (c_2)  \Gamma (a+b_2-c_2) \Gamma (-a+b_1-b_2+c_2)  }{\Gamma (a) \Gamma (b_1) \Gamma (b_2) \Gamma (-a-b_2+c_1+c_2)}\nonumber\\
		&\times\tilde{F}{}^{1:1;2}_{1:0;1}
  \left[
   \setlength{\arraycolsep}{0pt}
   \begin{array}{c@{{}:{}}c@{;{}}c}
  a+b_2-c_2 & c_1-b_1 & 1-b_2,c_2-b_2\\[1ex]
  a-b_1+b_2-c_2+1 & - & -a-b_2+c_1+c_2
   \end{array}
   \;\middle|\;
\frac{y-1}{x+y-1},\frac{x+y-1}{y}
 \right]\nonumber\\
		&+  (-x-y+1)^{-a}\left\langle\left(-\frac{y}{x+y-1}\right)^{b_2-c_2}\right\rangle\left\langle \left(\frac{y-1}{x+y-1}\right)^{-a+b_1-b_2+c_2}\right\rangle\frac{\Gamma (c_1) \Gamma (c_2) \Gamma (a-b_1+b_2-c_2) }{\Gamma (a) \Gamma (b_2) \Gamma (c_1-b_1)}\nonumber\\
		&\times{F}{}^{1:1;2}_{1:0;1}
  \left[
   \setlength{\arraycolsep}{0pt}
   \begin{array}{c@{{}:{}}c@{;{}}c}
  -a-b_2+c_1+c_2 & b_1 & 1-b_2,c_2-b_2\\[1ex]
  -a+b_1-b_2+c_2+1 & - & -a-b_2+c_1+c_2
   \end{array}
   \;\middle|\;
\frac{y-1}{x+y-1},\frac{y-1}{y}
 \right]
	\end{align}
	Region of convergence:  $|\frac{x+y-1}{y}| <1\land | \frac{y-1}{x+y-1}| <1$ (see Fig. \ref{ROCS6} \textit{left}).

	\subsection*{Series representation $S_{17}$}
	$S_{17}$ is Eq.(\ref{F2_ACMB1}).
	\begin{align}
		S_{17}&=(-y)^{-a} \frac{\Gamma (c_2) \Gamma (b_2-a)}{\Gamma (b_2) \Gamma (c_2-a)} F{}^{2:1;0}_{1:1;0}
  \left[
   \setlength{\arraycolsep}{0pt}
   \begin{array}{c@{{}:{}}c@{;{}}c}
 a,a-c_2+1 & b_1&-\\[1ex]
  a-b_2+1 & c_1 & \linefill
   \end{array}
   \;\middle|\;
 -\frac{x}{y},\frac{1}{y}
 \right]\nonumber\\
		&+(-y)^{-b_2} \frac{ \Gamma (c_2) \Gamma (a-b_2)}{\Gamma (a) \Gamma (c_2-b_2)}\,H_2\left(a-b_2,b_1,b_2,b_2-c_2+1;c_1;x,-\frac{1}{y}\right)
	\end{align}
	Region of convergence:  $| x| <1\land | -\frac{1}{y}| <1\land (| x| +1) | -\frac{1}{y}| <1$ (see Fig. \ref{ROCS6} \textit{middle}).

	\subsection*{Series representation $S_{18}$}
\begin{figure}[h]
		\centering
		\includegraphics[scale=.24]{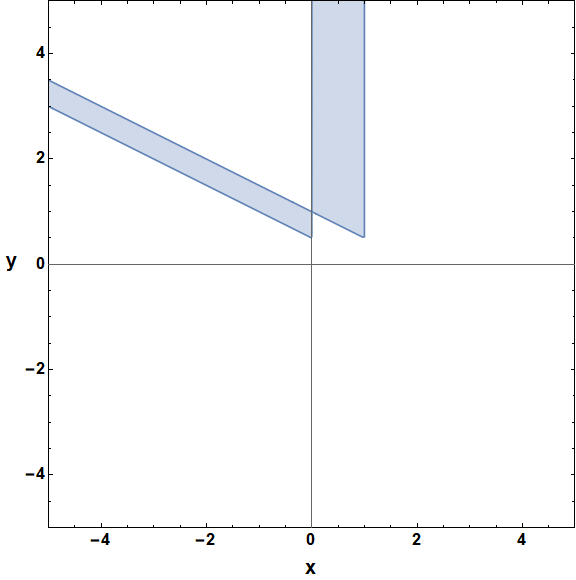}
		\includegraphics[scale=.24]{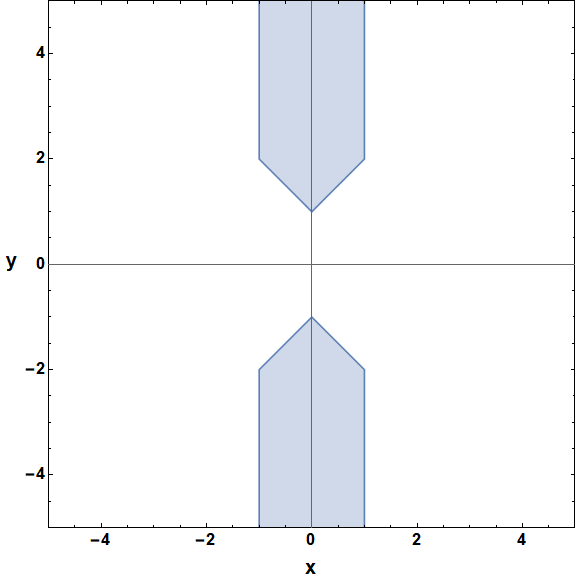}
		\includegraphics[scale=.24]{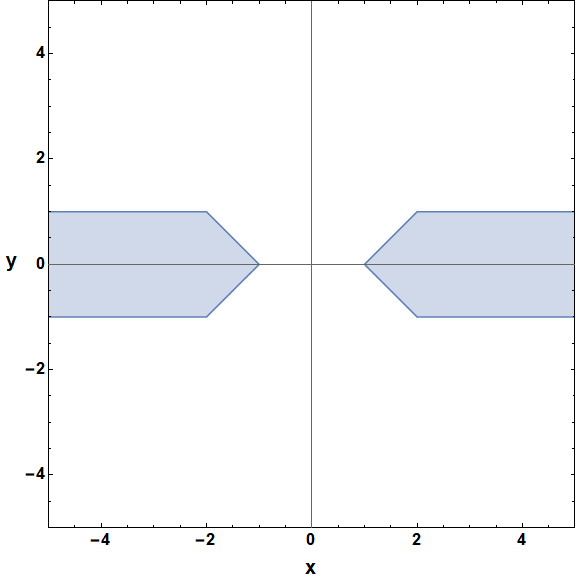}
		\caption{Convergence regions of $S_{16}$ (\textit{left}), $S_{17}$ (\textit{middle}) and $S_{18}$ (\textit{right}) for real values of $x$ and $y$.\label{ROCS6}}
	\end{figure}

	$S_{18}$  is the symmetrical partner of Eq.(\ref{F2_ACMB1}) obtained from Eq.(\ref{F_2sym}).
	\begin{align}
		S_{18}&=(-x)^{-a} \frac{\Gamma (c_1) \Gamma (b_1-a)}{\Gamma (b_1) \Gamma (c_1-a)} F{}^{2:1;0}_{1:1;0}
  \left[
   \setlength{\arraycolsep}{0pt}
   \begin{array}{c@{{}:{}}c@{;{}}c}
 a,a-c_1+1 & b_2&-\\[1ex]
  a-b_1+1 & c_2 & \linefill
   \end{array}
   \;\middle|\;
 -\frac{y}{x},\frac{1}{x}
 \right]\nonumber\\
		&+(-x)^{-b_1} \frac{ \Gamma (c_1) \Gamma (a-b_1)}{\Gamma (a) \Gamma (c_1-b_1)}\,H_2\left(a-b_1,b_2,b_1,b_1-c_1+1;c_2;y,-\frac{1}{x}\right)
	\end{align}
	Region of convergence:  $| y| <1\land \frac{1}{| x| }<1\land \frac{| y| +1}{| x| }<1$ (see Fig. \ref{ROCS6} \textit{right}).

\end{document}